\documentclass[11pt]{article}
\usepackage{amsmath,amssymb}
\usepackage{amsthm}
\usepackage{graphicx}
\usepackage{subfigure}
\usepackage{url}
\usepackage{multirow}
\usepackage[font=footnotesize]{subfig}
\numberwithin{equation}{section}

\usepackage{float}
\usepackage{array}
\usepackage{mathdots}
\usepackage{cellspace}
\usepackage[round]{natbib}
\usepackage[usenames,dvipsnames,svgnames,table]{xcolor}

\usepackage{a4wide}

\newtheorem{theorem2}{Theorem}[section]
\newtheorem{proposition}{Proposition}[section]
\newtheorem{lemma}{Lemma}[section]

\newtheorem{remark}{Remark}[section]

\newtheorem{assumption}{Assumption}

\begin{document}

\newcommand{\red}{\color{DarkRed}}
\newcommand\cov{\mathop{\text{cov}}}
\newcommand\var{\mathop{\text{var}}}
\newcommand\diag{\mathop{\text{diag}}}
\newcommand\cvd{\xrightarrow{~d~}{}}
\renewcommand{\(}{\left(}
\renewcommand{\)}{\right)}

\begin{center}
	{\bf \large On John's test for sphericity in large panel data models}
	\\[4mm]
	\title
	\author{Zhaoyuan LI\\[1mm] The Chinese University of Hong Kong, Shenzhen}\\[2mm] July 3, 2022
\end{center}

\bigskip

\begin{abstract}
	This paper studies John's test for sphericity of the error terms in large panel data models, where the number of cross-section units $n$ is large enough to be comparable to the number of times series observations $T$, or even larger. Based on recent random matrix theory results, John's test's asymptotic normality properties are established under both the null and the alternative hypotheses. These asymptotics are valid for general populations, i.e., not necessarily Gaussian provided certain finite moments. A fantastic phenomenon found in the paper is that John's test for panel data models possesses a powerful dimension-proof property. It keeps the same null distribution under different $(n,T)$-asymptotics, i.e., the small or medium panel regime $n/T\to 0$ as $T\to \infty$, the large panel regime $n/T\to c \in (0,\infty)$ as $ T\to \infty$, and the ultra-large panel regime $n/T\to \infty (T^\delta/n =O_p(1), 1<\delta<2)$ as $T\to \infty$. Moreover, John's test is always consistent except under the alternative of bounded-norm covariance with the large panel regime $n/T\to c \in (0,\infty)$.
	
\end{abstract}

\section{Introduction}
Consider the fixed-effects panel data regression model:
\begin{eqnarray}\label{model}
y_{it} = x_{it}^\prime \beta +\mu_i +\nu_{it}, \ \textrm{for} \ i=1,\ldots,n;\ t=1,\ldots,T,
\end{eqnarray}
where $i$ indexes the cross-sectional (individual) units and $t$ the time-series observations. The dependent variable is $y_{it}$ and $x_{it}$ denotes the exogenous regressors of dimension $k\times 1$, $\beta$ is the corresponding $k\times 1$ vector of parameters,  $\mu_i$
denotes the time-invariant individual effect, which could be correlated with the regressors $x_{it}$. Throughout the paper, the number of covariates $k$ is fixed while the two dimensions $n$ and $T$ may grow to infinity.
Let $\nu_t = (\nu_{1t},\ldots,\nu_{nt})^\prime$ be the panel-wise
error vector at time $t$. In this paper, we assume that the error vectors 
$\nu_1,\ldots, \nu_T$ are homoscedastic so that the cross-sectional covariance matrix
$\Sigma_n=\cov(\nu_t)$ is independent of $t$.

Our primary focus is on testing the sphericity of the covariance matrix of the disturbances:
\begin{eqnarray}\label{null}
H_0: \Sigma_n = \sigma_{\nu}^2 \mathbf{I}_n, \quad \textrm{vs} \quad H_1: \Sigma_n \neq \sigma_{\nu}^2 \mathbf{I}_n,
\end{eqnarray}
where $\sigma_{\nu}^2$ is an unknown positive constant and $\mathbf{I}_n$ the identity matrix of order $n$.

Testing for sphericity of symmetric errors in the linear model has a long history. Among others, the problems of testing for serial correlation, spatial correlation, unit roots, and stationarity can all be cast in this framework \citep{broda2019testing}. The null of sphericity expresses that the disturbances are cross-sectionally uncorrelated (independent if they are normal) and have the same variance (homoscedasticity). Rejecting the null means having cross-sectional correlation among the individual units of observations, heteroskedasticity, or both. This test is important for applied panel data work, given the numerous applications that report fixed effects estimation ignoring cross-section dependence or heteroscedasticity \citep{baltagi2011testing}. 

Much of the existing theory about the test for sphericity has been explored first in \cite{mauchly1940significance} about Gaussian likelihood ratio test and later in \cite{john1971some} and \cite{john1972distribution} about the invariant John's test, and also in textbooks like \cite{anderson2003introduction}. It is well known that classical multivariate procedures are generally challenged by large dimensional data \citep{paul2014random, yao2015sample}. The likelihood ratio test requires the dimension to be smaller than the sample size; otherwise, the likelihood ratio is identically null. Some corrections to traditional likelihood ratio test are proposed for large dimensional data; see \cite{wang2013sphericity} and \cite{li2016testing}. \cite{ledoit2002some} corrected John's test to cope with the high-dimensional context with normality assumption. \cite{wang2013sphericity} removed such a normality restriction and proved that the robustness of John's test is general. \cite{chen2010tests} proposed to use a family of well-selected U-statistics to test the sphericity; however, this test is very time-consuming \citep{li2016testing} and is illustrated to be slightly lower power than the corrected John's test proposed by \cite{wang2013sphericity}. \cite{li2016testing} studied John's test under an ultra-high-dimensional regime where the dimension is of a higher order than the sample size. 

These existing statistical tests for raw data do not apply directly to the sphericity test in panel regression models since the disturbances are unobservable. \cite{baltagi2011testing} extended \cite{ledoit2002some}'s result on John's test to the fixed effects panel data model. They corrected for the bias due to substituting the within residuals for the actual disturbances. \cite{baltagi2017asymptotic} studied the corresponding asymptotic power under the alternative of factor model. However, these results rely on the normality assumption of disturbances, a restriction for real data analysis . Building on the work of \cite{chen2010tests}, \cite{baltagi2015testing} extended the U-statistics to test for sphericity of the disturbances without the normality assumption. However, this kind of U-statistics carries a burden of extensive computations, and its consistency under the alternative is unknown.

This paper focuses on John's test due to its robustness proved by existing literature. We make the following main contributions.
\begin{enumerate}
	\item {\bf We consider general populations of disturbances}. Without normality restriction on disturbances, we derive the null distribution of John's test based on residuals for the sphericity test in (\ref{null}) and conduct power analysis. 
	
	\item {\bf We provide theoretical proof of the dimension-proof property of John's test in the large panel data model}. 
	This paper studies the panel data model (\ref{model}) in large dimensions, that is, the number of cross-section units is commensurate with or much larger than the number of times series observations. Accordingly, by imposing conditions on the relative speed at which $n$ and $T$ go to infinity, we define two asymptotic schemes:
	\begin{itemize}
		\item Large-panel asymptotic (LPA): the number of cross-sectional units $n$ and the number of time-series observations $T$ tend to infinity in such a way that $c_T = n/T \to c \in (0,\infty)$.
		\item Ultra-large panel asymptotic (ULPA): the number of cross-sectional units $n$ and the number of time-series observations $T$ tend to infinity in such a way that $c_T = n/T \to \infty$, $T^\delta /n = O_p(1), 1<\delta<2$.
	\end{itemize}
	The LPA scheme is appropriate for macroeconomic applications in which typically $n$ and $T$ are both large and of comparable magnitudes. The ULPA scheme is appropriate for microeconomic applications where typically $n$ is much larger than $T$. 
	
	John's test is studied under both of LPA and ULPA schemes defined above, and then its potent dimension-proof property is established, that is, John's test based on the residuals of the panel data model keeps the same null distribution under different $(n, T)$-asymptotics, i.e., the small or medium panel regime $n/T\to 0$ as $T\to \infty$, the large panel regime $n/T\to c \in (0,\infty)$ as $ T\to \infty$ and the ultra-large panel regime $n/T\to \infty, T^\delta/n =O_p(1), 1<\delta<2$ as $T\to \infty$.
	
	\item {\bf We conduct a comprehensive power analysis}. We study the asymptotic power of John's test under the factor model alternative with different settings. Without normality assumption on the disturbances and following \cite{baltagi2017asymptotic}, we investigate the asymptotic power of John's test for the panel model (\ref{model})  under the alternative of a weak factor model with fixed number of factors. Moreover, we conduct the power analysis under the following three scenarios, which are proposed as open questions by \cite{baltagi2017asymptotic}:
	\begin{enumerate}
			\item the asymptotic power as $\frac{n}{T}\to \infty$, i.e., under the ULPA scheme;
		\item the factors are neither strong nor weak;
		\item the factors are weak, but the number of factors goes to infinity jointly with $n$ and $T$. 
	\end{enumerate}
\end{enumerate}


The rest of the paper is organized as follows. Section 2 includes a literature review on John's test for raw data and panel regression model, respectively. The null distributions of residual-based John's test under both LAP and ULAP schemes for general populations of disturbances are presented in Section 3. Section 4 conducts the power analysis under the factor model alternative with different settings, followed by concluding remarks in Section 5. All technical proofs are relegated to the Appendix. 

\section{John's test, residual-based John's test}

\subsection{John's test for raw data}
\cite{john1971some} proposed a test statistic for the null hypothesis of sphericity in (\ref{null}) defined as follows:
\begin{eqnarray}\label{U}
U = \frac{1}{n} tr \left[\left(\frac{1}{n} tr S_ T\right)^{-1} S_T -\mathbf{I}_n \right]^2 = \left(\frac{1}{n} tr S_T \right)^{-2} \left(\frac{1}{n}tr S_T^2 \right) -1,
\end{eqnarray}
where $S_T$ is the {\em sample error covariance matrix} with order $n$. For the panel regression model (\ref{model}), $\displaystyle S_T=\frac{1}{T}\sum_{t=1}^T \nu_t \nu_t^\prime$. Notice that the statistic $U$ is independent of the scale parameter $\sigma_{\nu}^2$ under the null. Therefore, we assume w.l.o.g. $\sigma_{\nu}^2=1$ when dealing with the null distribution of $U$. 

\subsubsection{Large sample theory}

\cite{john1971some} proved that when the observations $\{\nu_{it}\}$ are normal, the sphericity test based on $U$ is a locally most powerful invariant test. Under these conditions, \cite{john1972distribution} established the large sample limiting distribution of $U$ under $H_0$, when $n$ is fixed and $T\to \infty$,
\begin{eqnarray}\label{john1972}
TU-n \overset{\mathcal{D}}{\longrightarrow} \frac{2}{n}\chi^2_{n(n+1)/2-1} -n.
\end{eqnarray}
This is referred to as John's test. It has been noticed for a while that John's test does not suffer from large dimensions, and this $\chi^2$ limit is quite accurate even when the ratio $n/T$ is not small.

\subsubsection{High-dimensional theory}
 \cite{ledoit2002some} studied the null distribution of $U$ under the LPA scheme. Assume that the observations $\{\nu_{it}\}$ are normal-distributed, under $H_0$ and the LPA scheme,
\begin{eqnarray}\label{ledoit}
TU-n \overset{\mathcal{D}}{\longrightarrow} N(1,4).
\end{eqnarray}
Meanwhile, if we let $n\to \infty$ on the right-hand side of (\ref{john1972}), it is not hard to see that 
\begin{eqnarray*}
	\frac{2}{n}\chi^2_{n(n+1)/2-1} -n \overset{\mathcal{D}}{\longrightarrow} N(1,4).
\end{eqnarray*}
In other words, \cite{ledoit2002some} extended the classical $T$-asymptotic theory (where $n$ is fixed) to the LPA scheme where $n$ goes to infinity proportionally with $T$.

Subsequently, \cite{wang2013sphericity} relaxed the normality restriction and proved that, assume $\{\nu_{it}\}$ are i.i.d., satisfying  $E \nu_{it}=0, E|\nu_{it}|^2=1, \gamma_4 = E|\nu_{it}|^4<\infty$, under $H_0$ and the LPA scheme, 
\begin{eqnarray}\label{Wangyao}
TU-n \overset{\mathcal{D}}{\longrightarrow} N(\gamma_4-2,4).
\end{eqnarray}

\subsubsection{Ultra high-dimensional theory}
\cite{li2016testing} further considered the asymptotic behavior of John's test statistic under the ULPA framework. Assume that $\{\nu_{it}\}$ are i.i.d., satisfying  $E (\nu_{it})=0, E|\nu_{it}|^2=1, \gamma_4 = E|\nu_{it}|^4<\infty$, under $H_0$ and the ULPA scheme,
\begin{eqnarray}\label{liyao}
TU-n \overset{\mathcal{D}}{\longrightarrow} N(\gamma_4-2,4),
\end{eqnarray}
which is the same as the result~(\ref{Wangyao}) under the LPA scheme. 

Since $\gamma_4=3+o(1)$ for normal distribution, the current results confirm each other. Besides, if the data has a non-normal distribution but has the same first four moments as the normal distribution, we have again $TU-n \overset{\mathcal{D}}{\longrightarrow} N(1,4)$ in (\ref{ledoit}). Another striking fact is that the limiting distribution of $U$ is independent of the dimension-to-sample ratio $\lim n/T$. John's test possesses a remarkable {\it dimension-proof} property, making it a very competitive candidate for sphericity testing regardless of the magnitude of $n,T$. Therefore, this paper focuses on John's test for sphericity test of the disturbances in the panel regression model (\ref{model}).

\subsection{John's test for panel data model}

In the panel regression model (\ref{model}), as $\{\nu_{it} \}$ are unobserved, the test statistic $U$ is modified using regression residuals. Consider the OLS estimator of the slope parameter $\beta$,
\begin{eqnarray*}
	\hat{\beta}= \left(\sum_{t=1}^T \sum_{i=1}^n \tilde{x}_{it}\tilde{x}_{it}^\prime \right)^{-1} \left(\sum_{t=1}^T \sum_{i=1}^n \tilde{x}_{it}\tilde{y}_{it} \right),
\end{eqnarray*}
where $\tilde{x}_{it} = x_{it}-\frac{1}{T}\sum_{t=1}^T x_{it}$ and $\tilde{y}_{it} = y_{it}-\frac{1}{T}\sum_{t=1}^T y_{it}$. It is well established that under mild conditions, $\hat{\beta}$ is a consistent estimator of $\beta$ no matter whether $n$ is fixed or tends to infinity jointly with $T$ \citep{baltagi2011testing,baltagi2017asymptotic}. The regression residuals $\{ \hat{\nu}_{it}\}$ are thus, 
\begin{eqnarray}
\hat{\nu}_{it} = \tilde{y}_{it} - \tilde{x}_{it}^\prime \hat{\beta} =\tilde{\nu}_{it} -\tilde{x}_{it}^\prime \left(\hat{\beta}-\beta \right),
\end{eqnarray}
where $\tilde{\nu}_{it}=\nu_{it}-\bar{\nu}_{i.}$ with $\bar{\nu}_{i.}=\frac{1}{T}\sum_{t=1}^T \nu_{it}$. The {\em sample residual covariance matrix} is $\displaystyle \hat{S}_T = \frac{1}{T}\sum_{t=1}^T \hat{\nu}_t \hat{\nu}_t^\prime$, where $\hat{\nu}_t = (\hat{\nu}_{1t},\ldots,\hat{\nu}_{nt})^\prime$ for $t=1,\ldots, T$, and the modified John's test statistic based on $\hat{S}_T$ is
\begin{eqnarray}\label{U_hat}
\hat{U} = \left(\frac{1}{n} tr \hat{S}_T \right)^{-2} \left(\frac{1}{n}tr \hat{S}_T^2 \right) -1.
\end{eqnarray}
which is referred to as the residual-based John's test statistic.

\subsubsection{The null distribution under LPA with normality}
 The asymptotics of $\hat{U}$ can be derived from the asymptotics of $U$ by estimating their difference, that is to examine the effect of replacing the errors $\{\nu_{it} \}$ with the residuals $\{\hat{\nu}_{it} \}$ on the asymptotics. Based on the result (\ref{ledoit}) of \cite{ledoit2002some}, \cite{baltagi2011testing} proved that, when the disturbances $\{\nu_{it} \}$ are normal-distributed, under $H_0$ and the LPA scheme, $\hat{U}$ follows asymptotically a normal distribution, 
\begin{eqnarray}\label{b2011}
T\hat{U}-n-\frac{n}{T-1} \overset{\mathcal{D}}{\longrightarrow} N(1,4).
\end{eqnarray}
Compared to (\ref{ledoit}), the use of residuals introduces a drift $\frac{n}{T-1}$ in the asymptotic mean of John's statistic. The test based on the asymptotic normal distribution given in the equation (\ref{b2011}) will be referred to as the {\em residual-based John's (RJ) test}.

\subsubsection{Asymptotic power of the RJ test}

\cite{baltagi2017asymptotic} studied the asymptotic power of the RJ test under the alternative of factor model, i.e., the disturbances $\{\nu_{it}\}$ follow a factor model, 
\begin{eqnarray}\label{factor_model}
\nu_{it} = \sum_{j=1}^r \xi_{ij}f_{tj} +\epsilon_{it},
\end{eqnarray}
where $\xi_{ij}$ is the factor loading of individual $i$ for factor $j$, $f_{tj}$ is the factor $j$ in period $t$, and $r$ is the known number of factors. To ensure the identification of the model, constraints are introduced as follows. Let $\xi_j = (\xi_{1j},\ldots,\xi_{nj})^\prime$, and $\Xi=(\xi_1,\ldots,\xi_r)$.
\begin{itemize}
	\item[(a)] The factors $\{f_{tj} \}$ are i.i.d. with mean $0$ and variance $\sigma_j^2<\infty$ across time, the correlation coefficient between factors $f_{tj}$ and $f_{tk}$ is zero, for all $j, k$ and $t$.
	\item[(b)] The vectors of factor loading $\xi_j$ are orthogonal to each other, that is, the matrix $\Xi^\prime \Xi$ is diagonal with distinct diagonal elements.
	\item[(c)] The idiosyncratic errors $\{\epsilon_{it} \}$ are i.i.d. $N(0, \sigma_{\epsilon}^2)$, and independent of all factors.
\end{itemize}
Then, the population covariance matrix of the disturbances $\{\nu_{it}\}$ is 
\begin{eqnarray} \label{Sigma_factor}
\Sigma_n = \sigma_\epsilon^2 \left(\sum_{j=1}^r h_j e_je_j^\prime +\mathbf{I}_n \right),
\end{eqnarray}
where $h_j = \frac{\sigma_j^2}{\sigma_\epsilon^2}||\xi_j||^2$, $e_j = \frac{\xi_j}{||\xi_j||}$, and $||e_j||=1$. And the covariance matrix can be viewed as a rank-$r$ perturbation of the null case. Let $\Sigma_n =\Gamma_n D_n \Gamma_n^\prime$, where $D_n = diag(\lambda_1, \ldots, \lambda_n)$ with the eigenvalues of $\Sigma_n$, $\lambda_j = \sigma_\epsilon^2 (1+h_j)$ for $j=1,\ldots, r$, and $\lambda_j=\sigma_\epsilon^2$ for $j=r+1,\ldots,n$. The matrix $\Gamma_n$ is the eigenvector matrix. This is known as spiked population model introduced in \cite{johnstone2001distribution}.

\cite{baltagi2017asymptotic} studied two kinds of alternatives: $h_j \to c_j \in (0, \infty)$ and $h_j/n \to c_j \in (0,\infty)$, which correspond to the weak and strong factor cases considered by \cite{onatski2012asymptotics} and  \cite{johnstone2009consistency}. They found that the RJ test is consistent under the strong factor alternative while  inconsistent under the weak factor alternative. Under the factor model alternative (\ref{factor_model}) with weak factors $h_j \to c_j \in (0, \infty)$, constraints (a)-(c), and the LPA scheme, 
\begin{eqnarray}\label{alter_dist_weak}
	T\hat{U}-n-\frac{n}{T-1}-\frac{T \sum_{j=1}^r h_j^2}{n} \overset{\mathcal{D}}{\longrightarrow} N(1,4).
\end{eqnarray}
Since $r$ is fixed, the limiting distribution of $\hat{U}$ under this alternative is just a rightward mean shift of its null distribution in (\ref{b2011}).

\section{Extension of the RJ test}
The extension is twofold: we relax the normality assumption on the disturbances and investigate the RJ test under both LPA and ULPA schemes. The powerful dimension-proof property of John's test for the panel data model is presented here. 

To derive the null distribution of the residual-based John's test statistic $\hat{U}$ in (\ref{U_hat}), we need the following assumptions to general populations of disturbances in the panel data model (\ref{model}).
\begin{assumption}\label{A1}
	The regressors $\{x_{it} \}$ are independent of the 
	disturbances $\{\nu_{it} \}$, and have uniformly bounded fourth moments, that is
	\[
	\sup_{i,t} E|x_{it}|^4\le C_1, \quad
	\text{for some positive constant } C_1.
	\]
\end{assumption}

\begin{assumption}\label{A2}
	The panel-wise error vectors $\nu_1,\ldots, \nu_T$ are i.i.d. with
	mean zero and  uniformly bounded eighth moments, that is
	\[
	\sup_{i,t} E|\nu_{it}|^8\le C_2, \quad
	\text{for some positive constant } C_2.
	\]
\end{assumption}

\begin{assumption}\label{A3}
	There is a doubly infinite array of i.i.d. random variables $(z_{ij}), i,j\geq 1$, satisfying
	\[E(z_{11})=0, \ E(z_{11}^2)=1, \ E(z_{11}^4) = \gamma_4 <\infty,\] 
	such that for each n, T, letting $Z_T=(z_{it} )_{1\leq i \leq n, 1\leq t \leq T}$, the idiosyncratic error vectors can be represented as $\nu_t = \Sigma_n^{1/2}z_j$ where $z_j=(z_{ij})_{1\leq i\leq n}$ denotes the $j$th column of $Z_T$.
\end{assumption}
\begin{assumption}\label{high_sigma}
	The population spectral distribution $H_n$ of $\Sigma_n$ weakly converges to a probability distribution $H$, as $n\to \infty$, and the sequence of spectral norms $(||\Sigma_n||)$ is bounded.
\end{assumption}
Assumption \ref{A1} is  standard in panel data models.
Under  Assumption \ref{A2}, the matrix $\Sigma_n=\cov(\nu_t)$ is constant over time. Under Assumptions \ref{A1} and \ref{A2}, $\hat{\beta}$ is a consistent estimator of $\beta$. Based on Assumption \ref{A3}, the sample error covariance matrix can be represented as $ S_T =T^{-1}\sum_{t=1}^T \nu_t \nu_t^\prime= T^{-1}\Sigma_n^{1/2}Z_TZ_T^\prime \Sigma_n^{1/2}$, which is motivated by the random matrix theory, and it is generic enough for a precise analysis of the sphericity test. Under the null, set $\sigma_{\nu}^2=1$, then $\Sigma_n=\mathbf{I}_n$ and thus $\nu_{it}=z_{it}$.

The following theorem establishes the asymptotic normality of $\hat{U}$ under the LPA scheme. \begin{theorem2}\label{theorem_high}
	Suppose Assumptions 1-4 hold for the panel data model (\ref{model}). Under $H_0$ and the LAP scheme, when $(n,T)\to \infty$, $c_T=n/T\to c \in(0, \infty),$  we have
	\begin{eqnarray}\label{null_gene}
	T\hat{U}-n - \left(\hat{\gamma}_4+c_T-2\right) \overset{\mathcal{D}}{\longrightarrow} N(0,4),
	\end{eqnarray}
where $\hat{\gamma}_4$ is the fourth moment of the residuals, i.e. $\hat{\gamma}_4=\frac{1}{nT}\sum_{i=1}^n \sum_{t=1}^T |\hat{\nu}_{it}|^4$, and $\hat{\gamma}_4=\gamma_4+o_p(1)$ as in the following Proposition \ref{diff_between_U}.
\end{theorem2}
Since $\hat{\gamma}_4=3+o_p(1)$ for normally distributed disturbances, this result reduces to the result (\ref{b2011}) for the normal case. In addition, if the errors have a non-normal distribution but have the same first fourth moments as the normal distribution, we have (\ref{b2011}) again as in the normal case. Therefore, Theorem \ref{theorem_high} constitutes a natural extension of the null distribution of $\hat{U}$ from Gaussian context to general population ones. The test based on the asymptotic normal distribution given in the equation (\ref{null_gene}) will be referred to as the {\em general residual-based John's (GRJ) test}.

Theorem \ref{theorem_high} is obtained by the combination of the result (\ref{Wangyao}), a general CLT for $U$ of the sample error covariance matrix of (unobserved) disturbances $\{\nu_{it}\}$, and the difference between $U$ and $\hat{U}$ of the sample residual covariance matrix from residuals $\{\hat{\nu}_{it} \}$ in the following proposition.
\begin{proposition}\label{diff_between_U}
	With the same assumptions as in Theorem \ref{theorem_high}, we have
	\begin{eqnarray}
	T(\hat{U}-U) - c_T =O_p\left(\frac{1}{n} \right),\label{prop1}\\
	\hat{\gamma_4}-\gamma_4 = O_p\left(\frac{1}{T} \right).\label{prop2}
	\end{eqnarray}
\end{proposition}
The proof of (\ref{prop1}) is similar to that of Proposition 4.3 in \cite{baltagi2011testing} for normal disturbances but is recalled in the Appendix together with the proof of (\ref{prop2}). This proposition implies that for the residual-based John's test statistic $\hat{U}$, the additional noise contained in the within residuals $\{\hat{\nu}_{it} \}$ will accumulate into a constant $c_T$, and the fourth moment $\hat{\gamma}_4$ of within residuals is a consistent estimator of the fourth moment $\gamma_4$ of the disturbances.

The following theorem establishes the asymptotic normality of $\hat{U}$ under the ULPA scheme.

\begin{theorem2}\label{ultra_null}
	Suppose Assumptions 1-4 hold for the panel data model (\ref{model}). Under $H_0$ and the ULPA scheme, when $(n, T)\to \infty, c_T=n/T\to \infty, T^\delta/n=O_p(1), 1<\delta<2$, we have
	\begin{eqnarray}
	T\hat{U}-n-(\hat{\gamma}_4+c_T-2) \overset{\mathcal{D}}{\longrightarrow} N(0,4).
	\end{eqnarray}
\end{theorem2}

Theorem \ref{ultra_null} is obtained by the combination of the result  (\ref{liyao}), a general CLT for $U$ of the sample error covariance matrix under ULPA, and the difference between $U$ and $\hat{U}$ of the sample residual covariance matrix. The proof is skipped as it is the same as that of Theorem \ref{null_gene}, i.e., the results under the LPA scheme in Proposition~\ref{diff_between_U} still hold under the ULPA scheme.

Theorem \ref{theorem_high} and Theorem \ref{ultra_null} together indicate the {\it dimension-proof} property of the GRJ test (or John's test for panel data model), i.e., regardless of normality, under any $(n,p)$-asymptotic, $n/T\to (0,\infty]$, has been testified. From the application perspective, where $n$ and $T$ are explicitly known, it becomes tough to decide which asymptotic scheme to use, namely, ``$n$ fixed, $T\to \infty$", ``$n/T\to c\in(0,\infty)$, $n,T\to \infty$", or ``$n/T\to \infty$, $n,T\to \infty$"etc. Therefore, the GRJ test possesses the powerful dimension-proof property, making it recommended over other procedures, as most test statistics constructed under the ``$n$ fixed, $T\to \infty$" asymptotic scheme would have different limiting distributions under ``$n/T\to c\in(0,\infty)$" and ``$n/T\to \infty$" asymptotic schemes.

\section{Asymptotic power}

In this section, we comprehensively study the GRJ test's asymptotic power. This is done under the alternative hypothesis of the factor model~(\ref{factor_model}) with the following three different settings:
\begin{itemize}
	\item[S1.] weak factors with a fixed number of factors;
	\item[S2.] weak factors with a divergent number of factors;
	\item[S3.] neither weak nor strong factors with a fixed number of factors.
\end{itemize}
The setting S1 has been studied by \cite{baltagi2017asymptotic} for normal-distributed disturbances. Without normality assumption, the constraints for the factor model (\ref{factor_model}) keep the same except that
\begin{itemize}
	\item[(c')] The idiosyncratic errors $\{\epsilon_{it} \}$ are i.i.d. $(0, \sigma_{\epsilon}^2)$, and independent of all factors.
\end{itemize}
 Under setting S1, the factors are weak and the number of factors is fixed, so the spectral norm of the covariance matrix~(\ref{Sigma_factor}) is bounded, while its spectral norm  is unbounded under the other two settings, S2 and S3. In the literature of random matrix theory, the asymptotic results of linear spectral statistics of covariance matrix with unbounded spectral norm are different from that of covariance matrix with bounded spectral norm. Therefore, the setting S1 and settings S2-S3 will be studied, respectively, in the following two subsections. 

\subsection{Covariance matrix with bounded spectral norm}
We first study the asymptotic power of the GRJ test under the alternative with setting S1 under the LPA scheme. 
\begin{theorem2}\label{high_weak}
	Suppose Assumptions 1, 3 and 4 hold for the panel data model (\ref{model}), and the disturbances follow the factor model (\ref{factor_model}) with constraints (a), (b) and (c'). Under the LAP scheme and the weak factor alternative S1 with $h_j\to d_j\in (0,\infty)$, $j=1,\ldots,r$, when $(n, T)\to \infty$, $c_T=n/T\to c\in (0, \infty)$, we have
	\begin{eqnarray}
		T\hat{U} -n - \left( \hat{\gamma}_4+c_T-2 \right) -  \frac{1}{c_T}\sum_{j=1}^r h_j^2 \overset{\mathcal{D}}{\longrightarrow} N(0,4).
	\end{eqnarray}
\end{theorem2}
Theorem \ref{high_weak} implies that under the weak factor alternative with a fixed number of factors, the limiting distribution of the residual-based John's test statistic is just a mean shift of its null distribution (\ref{null_gene}). When the disturbances are normal or $\gamma_4\to 3$, this result reduces to (\ref{alter_dist_weak}) for normal population.  From Theorem \ref{theorem_high} and Theorem \ref{high_weak}, it is easy to obtain the power function of the GRJ test under the LPA scheme, that is,
\begin{eqnarray}\label{eqn_power_weak}
	P_1 (H_1) = 1-\Phi \left(Z_\alpha -\frac{1}{2c_T} \sum_{j=1}^r h_j^2 \right),
\end{eqnarray}
which coincides with the result for the normal case. $P_1(H_1)\geq \alpha$, and the power of the test increases as spikes $\{h_j\}$ become large but decreases as the ratio $c_T$ of $n$ and $T$ increases.  Therefore, under the LAP scheme, the GRJ test is inconsistent under setting S1 regardless of the distribution of the disturbances.

Theorem \ref{high_weak} is obtained by combining the limiting distribution of $U$ and the difference between $U$ and $\hat{U}$ in the following propositions. 

\begin{proposition}\label{U_weak_alter}
	With the same assumptions as in Theorem \ref{high_weak}, the asymptotic normality of $U$ is
	\begin{eqnarray*}
		T{U} -n - \left( \frac{1}{c_T}\sum_{j=1}^r h_j^2 -{\gamma}_4+2 \right)  \overset{\mathcal{D}}{\longrightarrow} N(0,4).
	\end{eqnarray*}
\end{proposition}
\begin{proposition}\label{diff_between_U_weak}
	With the same assumptions as in Theorem \ref{high_weak}, we have
	\begin{eqnarray}\label{prop3}
		T(\hat{U}-U) - c_T =O_p\left(\frac{1}{n} \right),\\
		\hat{\gamma_4}-\gamma_4 = O_p\left(\frac{1}{T} \right).
	\end{eqnarray}
\end{proposition}

Next, we study the asymptotic power of the GRJ test under the alternative with setting S1 under the ULPA scheme.  

\begin{theorem2}\label{ultra_high_weak}
	Suppose Assumptions 1,3 and 4 hold for the panel data model (\ref{model}), and the disturbances follow the factor model (\ref{factor_model}) with constraints (a), (b) and (c'). In addition to Assumption \ref{high_sigma} of $\Sigma_n$, assume the following limits exits,
	\begin{itemize}
		\item[(a)] $\eta_1 = \lim_{n\to \infty} \frac{1}{n} tr (\Sigma_n)$,
		\item[(b)] $\eta_2= \lim_{n\to \infty} \frac{1}{n} tr(\Sigma_n^2)$,
		\item[(c)] $\eta_3 = \lim_{n\to \infty} \frac{1}{n}\sum_{i=1}^n (\Sigma_{ii})^2$, where $\Sigma_{ii}$ denotes the $i$-th diagonal entry of $\Sigma_n$.
	\end{itemize}
 Under the ULAP scheme and the weak factor alternative S1 with $h_j\to d_j\in (0,\infty)$, $j=1,\ldots,r$, when $(n, T)\to \infty$, we have
	\begin{eqnarray}
		T\hat{U} -n - \left(\frac{\eta_2}{\eta_1^2} -1\right) T - \frac{\eta_2+\eta_3(\hat{\gamma}_4-3)}{\eta_1^2} -c_T \overset{\mathcal{D}}{\longrightarrow} N\left(0,\frac{4\eta_2^2}{\eta_1^4}\right).
	\end{eqnarray}
\end{theorem2}
Note that if let $\Sigma_n = \sigma_{\nu}^2\mathbf{I}_n$, then $\eta_1=\sigma_{\nu}^2, \eta_2=\eta_3=\sigma_{\nu}^4$, this result reduces to the null distribution of John's test statistic under the ULPA scheme in Theorem \ref{ultra_null}.

From Theorem \ref{ultra_null} and Theorem \ref{ultra_high_weak}, it is easy to obtain the power function of the GRJ test under the ULPA scheme, that is,
	\begin{eqnarray*}
		P_2(H_1) \sim 1- \Phi \left(\frac{\eta_1^2}{\eta_2}Z_\alpha +\frac{\eta_1^2(\hat{\gamma}_4-2)-\eta_2-\eta_3(\hat{\gamma}_4-3)}{2\eta_2} +\frac{(\eta_1^2-\eta_2)T}{2\eta_2} \right) \to 1,
	\end{eqnarray*}
as $\eta_1^2<\eta_2$. Therefore, under setting S1, the GRJ test is consistent under the ULPA scheme.

Theorem \ref{ultra_high_weak} is obtained by combining the limiting distribution of $U$ and the difference between $U$ and $\hat{U}$ in the following propositions. 

\begin{proposition}\label{U_ulpa}
	
With the same assumptions as in Theorem \ref{ultra_high_weak}, we have
	\begin{eqnarray}
	TU-n -\left(\frac{\eta_2}{\eta_1^2} -1\right) T -\frac{\eta_2+\eta_3(\gamma_4-3)}{\eta_1^2} \overset{\mathcal{D}}{\longrightarrow} N\left( 0, \frac{4\eta_2^2}{\eta_1^4} \right).
	\end{eqnarray}
\end{proposition}

\begin{proposition}\label{conjection_general_ulpa}
	With the same assumptions as in Theorem \ref{ultra_high_weak}, when $(n,T)\to \infty, n/T=c_T\to c \in (0,\infty)$, we have
		\begin{eqnarray}
			T(\hat{U}-U) - c_T =O_p\left(\frac{1}{n} \right),\\
			\hat{\gamma_4}-\gamma_4 = O_p\left(\frac{1}{T} \right).
		\end{eqnarray}
	\end{proposition}

In conclusion, under Setting S1, the GRJ test is inconsistent under the LPA scheme where $n/T\to c\in(0, \infty)$, but consistent under the ULPA scheme where $n/T \to \infty$.

\begin{remark}
	Based on the above power analysis under the weak factor alternative with the LPA scheme, we can investigate the power behavior of the GRJ test under a general form of covariance matrix with bounded spectral norm, which covers the weak factor case as a special case; see the Supplementary Material.
\end{remark}

\subsection{Covaraince matrix with unbounded spectral norm}

In this subsection, we study the asymptotic power of the GRJ test under the other two settings S2 and S3, where the spectral norm of the covariance matrix $\Sigma_n$ of disturbances is unbounded. To achieve this goal, we consider John's test as a way to test the null hypothesis $H_0$  against an alternative hypothesis of the form:
\[H_1^\ast:	{ spec}(\Sigma_n)= \underbrace{\{\lambda_1,\ldots\lambda_{r_n}\}}_{G_1}\cup\underbrace{\{\lambda_{r_n+1},\ldots,\lambda_n\}}_{G_2}, \]
where the number of factors $r_n$ in $G_1$ can keep fixed (i.e., $r_n=r$) or diverges as the dimensions $n, T$ increase.
\begin{assumption}\label{rn}
	The eigenvalues in $G_1$ satisfy $\displaystyle \sup_i \max_{1\leq i \leq r_n} \lambda_i =O(n^\alpha)$ and the size of $
	G_1$ is $r_n/n\to \tau$ for some $0\leq \alpha< 1$ and $0\leq \tau <1$.
\end{assumption}
In Assumption \ref{rn}, the eigenvalues (factors) are weak (or bounded) when $\alpha=0$, and the eigenvalues (factors) are neither weak nor strong when $0<\alpha<1$. The number of factors is fixed when $\tau=0$, otherwise divergent.

We first give the asymptotic Gaussian of $U$ under the alternative $H_1^\ast$.
\begin{theorem2}\label{J_H1}
Suppose that Assumption \ref{A3} with $E(|z_{11}|^6)<\infty$ and Assumption \ref{rn} hold when $4\alpha+\tau\leq 1$ or Assumption \ref{A3} with $E(|z_{11}|^{16})<\infty$ and Assumption \ref{rn} hold when $4\alpha+\tau>1$, then under $H_1^\ast$ and the LPA scheme, when $c_T=n/T \to c \in (0, \infty)$ as $n, T\to \infty$, we have
	\begin{eqnarray}
		T\left( U- \mu \right)/\sigma \stackrel{D}{\rightarrow} N(0, 1),
	\end{eqnarray}
	where
	\[\mu= \frac{c_T \left\{(\gamma_4-3) tr(\Sigma_n \circ \Sigma_n)+(T+1)tr(\Sigma_n^2) \right\}}{(tr \Sigma_n)^2} +c_T-1, \]
	and
	\[\sigma^2 \sim \frac{1}{c_T^2}\left\{\frac{4 T^4 \theta_2^2}{\theta_1^6}\omega_{1}-\frac{4 T^4 \theta_2}{\theta_1^5}\omega_{2} +\frac{T^4}{\theta_1^4}\omega_{3} \right\}, \]
	with
	\begin{eqnarray*}
		\theta_1 &=& tr \Sigma_n,\\
		\theta_2 &=& T^{-1}\left\{(\gamma_4-3) tr(\Sigma_n \circ \Sigma_n) +(tr \Sigma_n)^2 +(T+1)tr(\Sigma_n^2) \right\},\\
		\omega_{1} &=& T^{-1}\left\{(\gamma_4-3) tr(\Sigma_n \circ \Sigma_n)+2(tr \Sigma_n)^2 \right\},\\
		\omega_{2} &=& T^{-2}\left\{4 tr \Sigma_n^2 tr \Sigma_n + 2 (\gamma_4-3) tr(\Sigma_n \circ \Sigma_n)tr \Sigma_n +2T(\gamma_4-3) tr(\Sigma_n \circ \Sigma_n^2)+4Ttr \Sigma_n^3  \right\},\\
		\omega_{3} &=& T^{-3} \Big\{8 tr \Sigma_n^2 (tr \Sigma_n)^2 + 4 (\gamma_4-3) (tr \Sigma_n)^2tr(\Sigma_n \circ \Sigma_n) + 16 T tr \Sigma_n tr \Sigma_n^3\\
		&& \quad \quad + 4 T (tr \Sigma_n^2)^2 + 8 T (\gamma_4-3) tr(\Sigma_n \circ \Sigma_n^2) tr \Sigma_n + 4 T^2 (\gamma_4-3) tr(\Sigma_n^2 \circ \Sigma_n^2) + 8 T^2 tr \Sigma_n^4  \Big\}.
	\end{eqnarray*}
\end{theorem2}
Theorem \ref{J_H1} implies that the limiting distribution of $U$ under $H_1^\ast$ relies on the magnitudes of the terms involving $\Sigma_n$. To illustrate the consistency of John's test, we analyze the power functions under a simplified factor model, where the covariance matrix $\Sigma_n$ is diagonal, 
\[\Sigma_n=diag(\lambda_1, \ldots, \lambda_{r_n}, \lambda_{r_n+1}, \ldots, \lambda_p). \]
\begin{itemize}
	\item Under setting S2: the eigenvalues in $G_1$ are weak (bounded), but the size of $G_1$ (the number of factors) is divergent, i.e., $\alpha=0, 0<\tau<1$, all trace terms involving $\Sigma_n$, e.g. $tr \Sigma_n$, $tr \Sigma_n^2$, $tr(\Sigma_n \circ \Sigma_n) (=tr \Sigma_n^2)$, $tr \Sigma_n^3$ and $tr \Sigma_n^4$, are an order of $n$. Then it is easy to verify that
	\[\mu = c_T-1+\frac{B_2}{B_1^2}\(1+\frac{1+\nu_4}{T} \), \ \textrm{and} \ \sigma^2=O(n), \]
	where $B_1=\frac{1}{n}tr\Sigma_n$ and $B_2=\frac{1}{n}tr\Sigma_n^2$. The power of John's test under this setting is
	\[P_1(H_1^\ast)=1-\Phi\(\frac{2Z_\alpha+(1+\nu_4)(1-\tau^{-1})+n(1-\frac{B_2}{B_1^2})}{\sigma}. \) \]
	As $B_2>B_1^2$ since all eigenvalues can't be equal, $1-\frac{B_2}{B_1^2}<0$, $n(1-\frac{B_2}{B_1^2})/\sigma\to -\infty$ as $n\to \infty$, thus the power function converges to 1. Therefore, John's test is consistent.

	\item Under setting S3: the eigenvalues in $G_1$ are divergent, but the size of $G_1$ is fixed, i.e., $\alpha>0, \tau=0$, set $r_n=r$, we have $tr \Sigma_n=O( r n^{\alpha})$, $tr \Sigma_n^2=tr(\Sigma_n \circ \Sigma_n)=O(r n^{2\alpha})$, $tr \Sigma_n^3=O(r n^{3\alpha})$ and $tr \Sigma_n^4=O(r n^{4\alpha})$. Then it is easy to verify that
	\[ \mu=c_1 n, \ \textrm{and} \ \sigma^2=c_2 n^3, \]
	where $c_1$ and $c_2$ are some positive constants. The power of John's test is
	\[P_2(H_1^\ast)=1-\Phi\left( \frac{2 Z_{\alpha}+n+(\gamma_4-2)-T\mu}{\sigma}\right). \]
	As $\frac{-T\mu}{\sigma}=-O(\sqrt{n})$, this power function converges to 1 as $n\to \infty$.
	Therefore, John's test is consistent.
\end{itemize}
Therefore, John's test for raw data is consistent under the alternative $H_1^\ast$. 

Next, we examine the consistency of the GRJ test. 

\begin{theorem2}\label{divergent_factors}
	Suppose Assumptions \ref{A1}, \ref{A3} and \ref{rn} hold for the panel data model \eqref{model} and the disturbances $\{\nu_{it}\}$ follow the factor model \eqref{factor_model} with constraints (a), (b) and (c'). Under the LPA scheme and the alternative $H_1^\ast$ with a divergent number of weak factors, the GRJ test is consistent.
\end{theorem2}
Theorem \ref{divergent_factors} is obtained by combining Theorem \ref{J_H1} and the following proposition.
\begin{proposition}\label{diff_between_U_diverge_factor}
	With the same assumptions as in Theorem \ref{divergent_factors},
		\begin{equation*}
		T(\hat{U}-U)=c_3c_T+c_4\tau+O_p\left(\frac{1}{\sqrt{n}}\right),
	\end{equation*}
  where $c_3$ and $c_4$ are some constants. 
\end{proposition}
Proposition \ref{diff_between_U_diverge_factor} implies that under the weak factor with divergent number alternative, additional noise obtained in $T\hat{U}-TU$ is $O_p(1)$. This magnitude is smaller than $O_p(T)$, the magnitude of $TU$. Thus, $T\hat{U}-TU$ is asymptotically dominated by $TU$, leading to the GRJ test's consistency.

\begin{theorem2}\label{not_weak}
		Suppose that for the panel data model \eqref{model}, Assumption \ref{A1}, Assumption \ref{A3} with $E(|z_{11}|^6)<\infty$ and Assumption \ref{rn} hold when $4\alpha+\tau\leq 1$, or Assumption \ref{A1}, Assumption \ref{A3} with $E(|z_{11}|^{16})<\infty$ and Assumption \ref{rn} hold when $4\alpha+\tau>1$, the disturbances $\{\nu_{it}\}$ follow the factor model \eqref{factor_model} with constraints (a), (b) and (c'). Under the LPA scheme and the alternative $H_1^\ast$ with a fixed number of neither strong nor weak factors, the GRJ test is consistent.
\end{theorem2}
Theorem \ref{not_weak} is obtained by combining Theorem \ref{J_H1} and the following proposition.
\begin{proposition}\label{diff_between_U_not_weak}
	With the same assumptions as in Theorem \ref{not_weak},
	\begin{equation*}
		TU=O_p\left(Tn^{2\alpha-1}\right)+O_p(T).
	\end{equation*}
	When $0<\alpha<\frac{1}{2}$,
	\begin{equation*}
		T(\hat{U}-U)=c_T+O_p\left(\frac{1}{\sqrt{T}}\right),
	\end{equation*}
	When $\frac{1}{2}\leq\alpha<1$,
	\begin{equation*}
		T(\hat{U}-U)=O_p\left(n^{2\alpha-1}\right).
	\end{equation*}
	
\end{proposition}

Proposition \ref{diff_between_U_not_weak} implies that under neither strong nor weak factor model, the magnitude of the additional noise contained in $T(\hat{U}-U)$ is smaller than $O_p(Tn^{2\alpha-1})+O_p(T)$, the magnitude of $TU$. Thus $T\hat{U}=TU+T(\hat{U}-U)$ is asymptotically dominated by $TU$, leading to the consistency of the GRJ test.

\begin{remark}
	In Theorem \ref{J_H1}, the limiting result still holds when $\alpha=1$, which corresponds to the strong factor case. Thus, together with the results under the strong factor alternative in Theorems 4 \& 5 of \cite{baltagi2017asymptotic}, we can conclude that the GRJ test is consistent under the strong factor alternative regardless of the distribution of disturbances.
\end{remark}

\begin{remark}
	In Theorem \ref{J_H1}, the factors and the number of factors can divergent simultaneously with $(n, T)$, i.e., $\alpha>0$ and $\tau>0$. We call it the fourth setting: neither weak nor strong factors with a divergent number of factors. We believe that the GRJ test is consistent under this setting, although the proof is omitted. 
\end{remark}

So far, we have conducted the power analysis of the GRJ test under the alternative of unbounded-norm covariance with the LPA scheme. The corresponding power analysis for the ULPA scheme is missing as there is no asymptotic result for linear spectral statistics of unbounded-norm covariance under the ULPA scheme. But based on the above power analyses, it is natural to conjecture that the GRJ test is still consistent under the ULPA scheme as it is consistent under both bounded-norm covariance with the scheme and unbounded-norm covariance with the LPA scheme.

\section{Conclusion}
This paper studies John's test for sphericity of the covariance matrix in a large panel data model with general populations of disturbances. The limiting distributions of the residual-based John's test are derived under both LPA and ULPA schemes. The results show that John's test possesses the powerful dimension-proof property with different $(n,T)$- asymptotic, making it recommended over other procedures. Moreover, the asymptotic power of John's test is conducted comprehensively under the factor model alternative with different settings. John's test is always consistent under the ULPA scheme, divergent number of weak factors and fixed number of non-weak factors, except under the weak factor alternative with a fixed number of factors under the LPA scheme. 

However, when the cross-sectional units $n$ increases sufficiently fast with time-series size $T$, say $n$ has an order of $T^2$ or larger, the limiting distribution of the residual-based John's test statistic $\hat{U}$ is unknown, which goes beyond the scope of the existing literature including this paper, and will require developing new mathematical techniques worthing much research effect in the future.

\section*{Supplementary Material}
The Supplementary Material contains the power analysis of the GRJ test under the general form of covariance matrix with bounded spectral norm under the LPA scheme.

\section*{Acknowledgments}
The author thanks Jianfeng Yao and Weiming Li for their helpful comments on the contents. Li's research was supported by NSFC grant 11901492.
\bibliographystyle{plainnat}
\bibliography{reference}

\appendix

\section{Preliminary knowledge of random matrix theory}

For nonnegative definite $n\times n$ Hermitian matrices $A_n$, the empirical spectral distribution (ESD) of $A_n$ generated by its eigenvalues  $a_1, \ldots, a_n$ is defined by 
\[F^{A_n} = \frac{1}{n}\sum_{i=1}^n \delta_{a_i},\] 
where $\delta_a$ denotes the Dirac mass at $a$. Assume that the ESD $H_n = F^{\Sigma_n}$ of $\Sigma_n$ converges weakly to a probability distribution $H$ on $[0,\infty)$. It is well-known that the (random) ESD $F^{A_n}$ of $A_n$ converges to a nonrandom limiting spectral distribution (LSD) $F^{c,H}$, called a generalized Marcenko-Paster law.
\[F^{A_n}(f) =\int f(x) d F^{A_n}(x)= n^{-1} \sum_{i=1}^n f(a_i)\] 
for suitable functions $f$, are of central importance in the proof here.

\section{Proof of Theorem \ref{theorem_high}}

\subsection*{Proof of Proposition \ref{diff_between_U}}

Recall $\tilde{y}_{it}=\tilde{x}_{it}^\prime \beta +\tilde{\nu}_{it}$ and 
\[\hat{\beta}-\beta=\left(\sum_{i=1}^n \sum_{t=1}^T \tilde{x}_{it}\tilde{x}_{it}^\prime \right)^{-1}\sum_{i=1}^n \sum_{t=1}^T \tilde{x}_{it}\tilde{\nu}_{it}, \]
Under Assumptions \ref{A1} and \ref{A2}, the within estimator $\hat{\beta}$ is $\sqrt{nT}$-consistent, that is,
\begin{eqnarray}\label{beta_asy}
	\hat{\beta}-\beta = O_p\left(\frac{1}{\sqrt{nT}} \right).
\end{eqnarray}
Here we also need the following fact.
\begin{lemma}\label{lemma_w}
	Define $w_{it}=\tilde{x}_{it}^\prime (\hat{\beta}-\beta)$ and $w_{t}=(w_{it} \cdots w_{nt})^\prime$. Based on model setting and Assumptions, we have $w_{it}=O_p\left(\frac{1}{\sqrt{nT}} \right)$ holds for all $i, t$.
\end{lemma}

The first equation \eqref{prop1} in Proposition~\ref{diff_between_U} basically follows the same proof line from Proposition 4.3 of \cite{baltagi2011testing}. However, for non-normal cases, $E(v_{it}^3)=0$ might not hold, which makes some difference. Now we verify that parts (a), (g), and (i) of Lemma A.3 of  \cite{baltagi2011testing}, which are the terms that contain the third order of error, still hold under non-normal assumption since the other parts of Lemma A.3 are the same for normal and non-normal distributions.

Following the settings of \cite{baltagi2011testing},
\begin{eqnarray*}
	&& A_0=\bar{v}\bar{v}^{'},\\
	&& A_1=\frac{1}{T}\sum_{t=1}^{T}{\tilde{x}_t(\tilde{\beta}-\beta)\tilde{v}_t^{'}},\\
	&& A_2=A_1^{'},\\
	&& A_3=\frac{1}{T}\sum_{t=1}^{T}{\tilde{x}_t(\tilde{\beta}-\beta)(\tilde{\beta}-\beta)^{'}\tilde{x}_t^{'}}.\\
\end{eqnarray*}
\begin{lemma}\label{lemma_w}
	Under Assumptions 1-3, (a)$\frac{1}{n}tr(SA_1)=O_p(\frac{1}{T^2})+O_p(\frac{1}{nT})+O_p(\frac{1}{T\sqrt{nT}})$;
	(g)$\frac{1}{n}tr(SA_0)=O_p(\frac{1}{T})+O_p(\frac{n}{T^2})$; 
	(i)$\frac{1}{n}tr(A_0 A_1)=O_p(\frac{1}{T^2})$.
\end{lemma}

\subsection*{Proof of part (a)}
The cubic terms involving $v_{it}^3$ are
\begin{eqnarray*}
	\frac{1}{nT^2}\sum_{t=1}^{T}\sum_{i=1}^{n}v_{it}^{3}\tilde{x}_{it}^{\prime}\left(\tilde{\beta}-\beta\right)+\frac{1}{nT^3}\sum_{t=1}^{T}\sum_{s=1}^{T}\sum_{i=1}^{n}v_{is}^{3}\tilde{x}_{it}^{\prime}\left(\tilde{\beta}-\beta\right)=O_p\left(\frac{1}{T\sqrt{nT}}\right),
\end{eqnarray*}
which is smaller than the order of the leading terms of $\frac{1}{n}tr(SA_1)$, so $\frac{1}{n}tr(SA_1)=O_p(\frac{1}{T^2})+O_p(\frac{1}{nT})+O_p(\frac{1}{T\sqrt{nT}})$ still holds under non-normal assumption.

\subsection*{Proof of part (g)}
The cubic terms involving $v_{it}^3$ are
\begin{eqnarray*}
	\frac{1}{nT^3}\sum_{\tau\neq s}^{T}\sum_{s=1}^{T}\sum_{i=1}^{n}v_{is}^{3}v_{i\tau}+\frac{1}{nT^3}\sum_{t\neq s}^{T}\sum_{s=1}^{T}\sum_{i=1}^{n}v_{it}^{3}v_{is}=O_p\left(\frac{1}{T\sqrt{nT}}\right),
\end{eqnarray*}
which is smaller than the order of the leading terms of $\frac{1}{n}tr(SA_0)$, so $\frac{1}{n}tr(SA_0)=O_p(\frac{1}{T})+O_p(\frac{n}{T^2})$ still holds under non-normal assumption.

\subsection*{Proof of part (i)}
The cubic terms involving $v_{it}^3$ are 
\begin{eqnarray*}
	&& \frac{1}{nT^3}\sum_{s=1}^{T}\sum_{i=1}^{n}v_{is}^3\tilde{x}_{is}^{\prime}(\hat{\beta}-\beta)
	+\frac{1}{nT^4}\sum_{s=1}^{T}\sum_{\tau=1}^{T}\sum_{i=1}^{n}v_{is}^3\tilde{x}_{i\tau}^{\prime}(\hat{\beta}-\beta)\\
	&=& O_p(\frac{1}{T^2\sqrt{nT}})+O_p(\frac{1}{T^2\sqrt{nT}}),
\end{eqnarray*}
which is smaller than the oder of the leading term of $\frac{1}{n}tr(A_0 A_1)$, so $\frac{1}{n}tr(A_0 A_1)=O_p(\frac{1}{T^2})$ still holds under non-normal assumption.

For the first equation in Proposition \ref{diff_between_U}, following (A.8) of \cite{baltagi2011testing}, we have
\begin{eqnarray*}
	T(\hat{U}-U) &=& \frac{\frac{n}{T} -\frac{n}{T^2}+O_p\left(\frac{\sqrt{n}}{T} \right) + O_p\left(\frac{1}{\sqrt{n}} \right) +O_p\left(\frac{1}{\sqrt{T}} \right)}{\frac{(T-1)^2}{T^2} +O_p\left(\frac{1}{\sqrt{nT}} \right)} \\
	&=& \frac{Tn}{(T-1)^2} -\frac{n}{(T-1)^2}+O_p\left(\frac{\sqrt{n}}{T-1} \right) + O_p\left(\frac{1}{\sqrt{n}} \right) +O_p\left(\frac{1}{\sqrt{T}} \right),
\end{eqnarray*}
Obviously, $T(\hat{U}-U)-c_T\to 0$ as $(n,T)\to \infty$ with $c_T=n/T\to c \in (0, \infty)$.

For the second equation in Proposition \ref{diff_between_U}, as $\hat{\nu}_{it} = \nu_{it}-\bar{\nu}_{i.}-\tilde{x}_{it}^\prime \left(\hat{\beta}-\beta \right)$, we have
\begin{eqnarray}\label{gamma_diff}
	E|\hat{\nu}_{it}|^4 &=& E\left|\nu_{it} - \bar{\nu}_{i.}-\tilde{x}_{it}^\prime (\hat{\beta}-\beta)\right|^4  \nonumber\\
	&=& E(\nu_{it}-\bar{\nu}_{i.})^4 + E\left[\tilde{x}_{it}^\prime (\hat{\beta}-\beta)\right]^4 + 6 E (\nu_{it}-\bar{\nu}_{i.})^2\left[\tilde{x}_{it}^\prime (\hat{\beta}-\beta)\right]^2.
\end{eqnarray}
The first term of the sum in (\ref{gamma_diff}) is 
\begin{eqnarray*}
	E(\nu_{it}-\bar{\nu}_{i.})^4 &=& E \left(\nu_{it}^4 + \bar{\nu}_{i.}^4 -4 \nu_{it}^3\bar{\nu}_{i.} + 6 \nu_{it}^2\bar{\nu}_{i.}^2 -4\nu_{it}\bar{\nu}_{i.}^3  \right) \\
	&=& E(\nu_{it}^4)\left(1-\frac{4}{T}+\frac{6}{T^2}-\frac{3}{T^3} \right) + \frac{6(T-1)}{T^2} -\frac{11(T-1)}{T^3} \\
	&=& E(\nu_{it}^4) +O_p\left(\frac{1}{T}\right).
\end{eqnarray*}
The second term of the sum in (\ref{gamma_diff}) is
\begin{eqnarray*}
	E\left[\tilde{x}_{it}^\prime (\hat{\beta}-\beta)\right]^4 = O_p\left(\frac{1}{n^2T^2} \right),
\end{eqnarray*}
due to $\tilde{x}_{it}^\prime (\hat{\beta}-\beta) = O_p\left(\frac{1}{\sqrt{nT}} \right)$. And the third term of the sum in (\ref{gamma_diff}) is
\begin{eqnarray*}
	E (\nu_{it}-\bar{\nu}_{i.})^2\left[\tilde{x}_{it}^\prime (\hat{\beta}-\beta)\right]^2 = O_p\left(\frac{1}{nT} \right).
\end{eqnarray*}
Therefore,
\begin{eqnarray*}
	E|\hat{\nu}_{it}|^4 = E(\nu_{it}^4) +O_p\left(\frac{1}{T} \right).
\end{eqnarray*}
The proof of Proposition~\ref{diff_between_U} is complete. 

\subsection*{Proof of Theorem \ref{theorem_high}}
Following the result (\ref{Wangyao}) of \cite{wang2013sphericity}, John's test statistic $U$ is asymptotically normal,
\begin{eqnarray}\label{wangyao2}
TU-n -(\gamma_4-2) \overset{\mathcal{D}}{\longrightarrow} N(0,4).
\end{eqnarray}
Let $J = TU-n-(\gamma_4-2)$. 
Based on residuals, the statistic will be $\hat{J}=T\hat{U}-n-(\hat{\gamma}_4-2)$. Notice that 
\begin{eqnarray*}
	\hat{J} = J + (\hat{J}-J) = J + T(\hat{U}-U)-(\hat{\gamma_4}-\gamma_4).
\end{eqnarray*}
Consequently, the effect of replacing the error $\{\nu_{it} \}$ with the residuals $\{\hat{\nu}_{it} \}$ on the asymptotic is captured by the difference $\hat{J}-J=T(\hat{U}-U)-(\hat{\gamma}_4-\gamma_4)$. Based on Proposition~\ref{diff_between_U}, the proof of Theorem~\ref{theorem_high} is complete. 

\section{Proof of Theorem \ref{high_weak}}

\subsection*{Proof of Proposition \ref{U_weak_alter}}
Let $\{\ell_i \}_{1\leq i\leq n}$ be the eigenvalues of the sample error covariance matrix $S_T$. We first derive the limiting distribution of $\left(\frac{1}{n}tr S_T^2,\frac{1}{n}tr S_T \right)$. 
Let $f(x)=x^2$ and $g(x)=x$. From the decompositions
\begin{eqnarray*}
	tr S_T^2=\sum_{i=1}^n \ell_i^2 &=& n \int f(x) d\left(F^{S_T}(x)-F^{c_T,H_n}(x) \right) + n  F^{c_T,H_n}(f),\\
 tr S_T=\sum_{i=1}^n \ell_i &=& n \int g(x)d\left(F^{S_T}(x)-F^{c_T,H_n}(x) \right) + n  F^{c_T,H_n}(g),
\end{eqnarray*}
define 
\begin{eqnarray*}
	X_n(f) &=& \sum_{i=1}^n \ell_i^2 - n  F^{c_T,H_n}(f),\\
	X_n(g) &=& \sum_{i=1}^n \ell_i - n  F^{c_T,H_n}(g).
\end{eqnarray*}
According to the CLT for linear spectral statistics established in \cite{pan2008central} for general populations, as $n/T=c_T\to c\in(0,\infty)$, we have
\begin{eqnarray}\label{CLT_Xn}
\left(\begin{array}{c}
X_n(f) \\ X_n(g) 
\end{array} \right)
\overset{\mathcal{D}}{\longrightarrow} N \left( \boldsymbol{\mu}, \mathbf{V} \right),
\end{eqnarray}
with
\begin{eqnarray*}
	\boldsymbol{\mu} = \left(\begin{array}{c}
		(\gamma_4-2)c \\ 0
	\end{array} \right),
\end{eqnarray*}
and
\begin{eqnarray*}
	\mathbf{V} = \left(\begin{array}{cc}
		4c^2+4(\gamma_4-1)(c+2c^2+c^3) & 2(\gamma_4-1)(c+c^2) \\
		2(\gamma_4-1)(c+c^2) & (\gamma_4-1)c
	\end{array} \right)
\end{eqnarray*}
The spikes (factors) do not prevent the above central limit theorem. However, the centering term will be modified according to the values of the spikes. From Theorem 2 in \cite{wang2014note}, we obtain the centering terms in $X_n(f)$ and $X_n(g)$, respectively, 
\begin{eqnarray*}
	F^{c_T,H_n}(f) &=& 1+c_T+ \frac{2}{T} \sum_{j=1}^r (1+h_j) + \frac{1}{n}\sum_{j=1}^r (1+h_j)^2  -\frac{2r}{T}-\frac{r}{n} + O\left(\frac{1}{T^2} \right),\\
	F^{c_T,H_n}(g) &=&1-\frac{r}{n}+ \frac{1}{n}\sum_{j=1}^r (1+h_j) +O\left(\frac{1}{T^2} \right).
\end{eqnarray*}
Therefore, the result (\ref{CLT_Xn}) can be rewritten as
\begin{eqnarray*}
	T \left(\begin{array}{c} n^{-1} \sum_{i=1}^n\ell_i^2 -m_1 \\ n^{-1}\sum_{i=1}^n \ell_i -m_2 \end{array} \right) \overset{\mathcal{D}}{\longrightarrow} N\left(\left(\begin{array}{c}
		0 \\ 0
	\end{array} \right), \frac{1}{c^2}\cdot \mathbf{V} \right),
\end{eqnarray*}
with 
\begin{eqnarray*}
	m_1 &=& 1+c_n+ \frac{2}{T} \sum_{j=1}^r (1+h_j) + \frac{1}{n}\sum_{j=1}^r (1+h_j)^2  -\frac{2r}{T}-\frac{r}{n} -\frac{(\gamma_4-2)c}{n}, \\
	m_2 &=& 1-\frac{r}{n}+ \frac{1}{n}\sum_{j=1}^r (1+h_j).
\end{eqnarray*}
Define the function $q(x,y)=\frac{x}{y^2}-1$, then $U = q\left( n^{-1} trS_T^2, n^{-1} tr S_T \right)$. We have
\begin{eqnarray*}
	&& \frac{\partial q}{\partial x}(m_1, m_2)=m_2^{-2},\\
	&& \frac{\partial q}{\partial y}(m_1, m_2) = -2m_1 m_2^{-3},\\
	&& q(m_1,m_2) = m_1 m_2^{-2}-1.
\end{eqnarray*}
By the delta method,
\begin{eqnarray*}
	T (U - q(m_1, m_2)) \overset{\mathcal{D}}{\longrightarrow} N(0, \lim V_U),
\end{eqnarray*}
where
\begin{eqnarray*}
	V_U = \left(\begin{array}{c} \frac{\partial q}{\partial x}(m_1, m_2) \\ \frac{\partial q}{\partial y}(m_1, m_2) \end{array} \right)^T \cdot \frac{1}{c^2}\mathbf{V} \cdot \left(\begin{array}{c} \frac{\partial q}{\partial x}(m_1, m_2) \\ \frac{\partial q}{\partial y}(m_1, m_2) \end{array} \right) \to 4.
\end{eqnarray*}
Therefore, 
\begin{eqnarray*}
	T \left(U- q(m_1,m_2) \right) \overset{\mathcal{D}}{\longrightarrow} N(0,4),
\end{eqnarray*}
where 
\begin{eqnarray*}
	q(m_1,m_2) = c_T +\frac{1}{n}\sum_{j=1}^r h_j^2-\frac{(\gamma_4-2)c}{n}.
\end{eqnarray*}
The proof of Proposition \ref{U_weak_alter} is complete.

\subsection*{Proof of Proposition \ref{diff_between_U_weak}} 
\begin{lemma}\label{lemma_beta}
	Under Assumptions 1, 2 and the weak factor alternative with $h_j\rightarrow d_j\in(0,\infty)$ for $j = 1,\ldots,r$
	\begin{equation*}
		\tilde{\beta}-\beta=O_p(\frac{1}{\sqrt{nT}})
	\end{equation*}
\end{lemma}

\subsubsection*{Proof of Lemma \ref{lemma_beta}}
Recall $\tilde{y}_{it}=\tilde{x}_{it}^{\prime}\beta+\tilde{v}_{it}$ and
\begin{equation*}
	\tilde{\beta}-\beta=\left(\sum_{t=1}^{T}\sum_{i=1}^{n}\tilde{x}_{it}\tilde{x}_{it}^{\prime}\right)^{-1} \sum_{t=1}^{T}\sum_{i=1}^{n}\tilde{x}_{it}\tilde{v}_{it},
\end{equation*}
where $\left(\sum_{t=1}^{T}\sum_{i=1}^{n}\tilde{x}_{it}\tilde{x}_{it}^{\prime}\right)^{-1}=O_p(\frac{1}{nT})$. The latter part can be written as
\begin{equation}\label{c2}
	\begin{aligned}
		\sum_{t=1}^{T}\sum_{i=1}^{n}\tilde{x}_{it}\tilde{v}_{it}
		=\sum_{t=1}^{T}\sum_{i=1}^{n}\tilde{x}_{it}v_{it}-\frac{1}{T}\sum_{t=1}^{T}\sum_{i=1}^{n}\sum_{s=1}^{T}\tilde{x}_{it}{v}_{is}.\\
	\end{aligned}
\end{equation}
The first term in (\ref{c2}) is $\sum_{t=1}^{T}\sum_{i=1}^{n}\tilde{x}_{it}v_{it}=O_p(\frac{1}{\sqrt{nT}})$ because
\begin{equation*}
	\begin{aligned}
		E\left(\sum_{t=1}^{T}\sum_{i=1}^{n}\tilde{x}_{it}v_{it}\right)^2 &= \sum_{j=1}^{n}\sum_{t=1}^{T}\sum_{i=1}^{n}E\tilde{x}_{it}v_{it}\tilde{x}_{jt}v_{jt}+\sum_{s\neq t}^{T}\sum_{j=1}^{n}\sum_{t=1}^{T}\sum_{i=1}^{n}E\tilde{x}_{it}v_{it}\tilde{x}_{js}v_{js}\\
		&=\sum_{j=1}^{n}\sum_{t=1}^{T}\sum_{i=1}^{n}E\tilde{x}_{it}\tilde{x}_{jt}E v_{it}v_{jt}+0\\
		&=\sum_{t=1}^{T}\sum_{i=1}^{n}E\tilde{x}_{it}^2E v_{it}^2+\sum_{j\neq i}^{n}\sum_{t=1}^{T}\sum_{i=1}^{n}E\tilde{x}_{it}\tilde{x}_{jt}E v_{it}v_{jt}\\
		&=O(nT)+\sum_{j\neq i}^{n}\sum_{t=1}^{T}\sum_{i=1}^{n}E\tilde{x}_{it}\tilde{x}_{jt}E v_{it}v_{jt}\\
		&\leq O(nT)+\sqrt{C_1}\sum_{j\neq i}^{n}\sum_{t=1}^{T}\sum_{i=1}^{n}E v_{it}v_{jt}\\
		&= O(nT)+\sqrt{C_1}\sum_{j\neq i}^{n}\sum_{t=1}^{T}\sum_{i=1}^{n}\sigma^2\left(\sum_{k=1}^{r}h_k e_{ik} e_{jk}\right)\\
		&\leq O(nT)+T\sqrt{C_1}\sigma^2\left(\sum_{k=1}^{r}h_k\sum_{i=1}^{n}e_{ik}\sum_{j=1}^{n}e_{jk}\right)\\
		&\leq O(nT)+T\sqrt{C_1}\sigma^2\left(\sum_{k=1}^{r}h_k\sqrt{n\sum_{i=1}^{n}e_{ik}^2}\sqrt{n\sum_{j=1}^{n}e_{jk}}\right)\\
		&\leq O(nT)+T\sqrt{C_1}\sigma^2\left(h\sqrt{n}\sqrt{n}\right)\\
		&=O(nT).
	\end{aligned}
\end{equation*}

In the calculation above, w.l.o.g, we treat $\tilde{x}_{it}$ as a scalar. We apply inequality $\left(\sum_{i=1}^{n}a_i\right)^2\leq n\sum_{i=1}^{n}a_i^2$ , $E\tilde{x}_{it}\tilde{x}_{jt}\leq \sup_{i,t}E_{x_{it}^2}\leq\sqrt{E x_{it}^4}\leq \sqrt{C_1}$ and $\Vert e_k\Vert=\sum_{i=1}^{n}e_{ik}^2=1$.

The second term in (\ref{c2}) is $\frac{1}{T}\sum_{t=1}^{T}\sum_{i=1}^{n}\sum_{s=1}^{T}\tilde{x}_{it}{v}_{is}=O_p(\sqrt{nT})$ because

\begin{equation*}
	\begin{aligned}
		E\left(\frac{1}{T}\sum_{t=1}^{T}\sum_{i=1}^{n}\sum_{s=1}^{T}\tilde{x}_{it}{v}_{is}\right)^2
		&=\frac{1}{T^2}\sum_{t_2=1}^{T}\sum_{i_2=1}^{n}\sum_{s_2=1}^{T}\sum_{t_1=1}^{T}\sum_{i_1=1}^{n}\sum_{s_1=1}^{T}
		E\tilde{x}_{i_1 t_1}\tilde{x}_{i_2 t_2}E{v}_{i_1 s_1}{v}_{i_2 s_2}\\
		&=\frac{1}{T^2}\sum_{t_2=1}^{T}\sum_{i_2=1}^{n}\sum_{t_1=1}^{T}\sum_{i_1=1}^{n}\sum_{s_1=1}^{T}
		E\tilde{x}_{i_1 t_1}\tilde{x}_{i_2 t_2}E{v}_{i_1 s_1}{v}_{i_2 s_1}\\
		&\leq \frac{\sqrt{C}}{T^2}\sum_{t_2=1}^{T}\sum_{i_2=1}^{n}\sum_{t_1=1}^{T}\sum_{i_1=1}^{n}\sum_{s_1=1}^{T}
		E{v}_{i_1 s_1}{v}_{i_2 s_1}\\
		&=\sqrt{C}\sum_{i_2=1}^{n}\sum_{i_1=1}^{n}\sum_{s_1=1}^{T}E{v}_{i_1 s_1}{v}_{i_2 s_1}\\
		&= O(nT),
	\end{aligned} 
\end{equation*}
where $C$ is a constant. Thus $\sum_{t=1}^{T}\sum_{i=1}^{n}\tilde{x}_{it}\tilde{v}_{it}=O_p(\sqrt{nT})$. The proof of Lemma \ref{lemma_beta} is complete.

\begin{lemma}\label{diff_between_S}
	Under Assumptions 1, 2 and under the weak factor alternative with $h_j\rightarrow d_j\in(0,\infty)$ for $j = 1,\ldots,r$,
	
	(a)$\frac{1}{n}tr(S)=\sigma^2+O_p(\frac{1}{\sqrt{nT}})$;
	
	(b)$\frac{1}{n}tr(S^2)=(\frac{n}{T}+1)\sigma^4+O_p(\frac{1}{\sqrt{T}})$;
	
	(c) $\frac{1}{n}tr(\hat{S})-\frac{1}{n}tr(S)=-\frac{\sigma^2}{T}+O_p(\frac{1}{T\sqrt{n}})$;
	
	(d) $\frac{1}{n}tr(\hat{S}^2)-\frac{1}{n}tr(S^2)=-\frac{2}{T}\sigma^4-\frac{n}{T^2}\sigma^4+O_p(\frac{1}{T\sqrt{T}})$.
\end{lemma}

\subsubsection*{Proof of Lemma \ref{diff_between_S}}
Proof of part (a).
\begin{eqnarray*}
	\frac{1}{n}tr(S)
	&=&\frac{1}{nT}\sum_{t=1}^T\sum_{i=1}^{n}v_{it}^2=\frac{1}{nT}\sum_{t=1}^T\sum_{i=1}^{n}\left(v_{it}^2-Ev_{it}^2\right)+\frac{1}{n}\sum_{i=1}^n Ev_{it}^2\\
	&=&\frac{1}{nT}\sum_{t=1}^T\sum_{i=1}^{n}\left(v_{it}^2-Ev_{it}^2\right)+\frac{1}{n}\sum_{i=1}^n \sigma^2\left(1+\sum_{j=1}^r h_j e_{ij}^2\right)\\
	&=&\frac{1}{nT}\sum_{t=1}^T\sum_{i=1}^{n}\left(v_{it}^2-Ev_{it}^2\right)+\sigma^2+\frac{\sum_{j=1}^r h_j}{n}\\
	&=&\sigma^2+O_p(\frac{1}{\sqrt{nT}}),
\end{eqnarray*}
due to $\frac{1}{nT}\sum_{t=1}^T\sum_{i=1}^{n}\left(v_{it}^2-Ev_{it}^2\right)=O_p(\frac{1}{\sqrt{nT}})$, which is obtained from
\begin{eqnarray*}
	&&E\left(\frac{1}{nT}\sum_{t=1}^T\sum_{i=1}^{n}\left(v_{it}^2-Ev_{it}^2\right)\right)^2\\
	&=&\frac{1}{n^2 T^2}\sum_{t=1}^T\sum_{i=1}^{n}E\left(v_{it}^2-Ev_{it}^2\right)^2+\frac{1}{n^2 T^2}\sum_{j\neq i}^{n}\sum_{t=1}^T\sum_{i=1}^{n}E\left(v_{it}^2-Ev_{it}^2\right)\left(v_{jt}^2-Ev_{jt}^2\right)\\
	&=& O(\frac{1}{nT})+\frac{1}{n^2 T^2}\sum_{j\neq i}^{n}\sum_{t=1}^T\sum_{i=1}^{n}E\left(v_{it}^2-Ev_{it}^2\right)\left(v_{jt}^2-Ev_{jt}^2\right)\\
	&=& O(\frac{1}{nT})+\frac{1}{n^2 T^2}\sum_{j\neq i}^{n}\sum_{t=1}^T\sum_{i=1}^{n}E\left(v_{it}^2 v_{jt}^2\right)-Ev_{it}^2 Ev_{jt}^2\\
	&=& O(\frac{1}{nT})+\frac{1}{n^2 T^2}\sum_{j\neq i}^{n}\sum_{t=1}^T\sum_{i=1}^{n}
	\left(E(\gamma_i f_t+\epsilon_{it})^2(\gamma_j f_t+\epsilon_{jt})^2-(\gamma_i^2\sigma_1^2+\sigma^2)(\gamma_j^2\sigma_1^2+\sigma^2)\right)\\
	&=& O(\frac{1}{nT})+\frac{1}{n^2 T^2}\sum_{j\neq i}^{n}\sum_{t=1}^T\sum_{i=1}^{n}
	\left(E\left(\gamma_i^2 f_t^2+2\gamma_i f_t\epsilon_{it}+\epsilon_{it}^2\right)\left(\gamma_j^2 f_t^2+2\gamma_j f_t\epsilon_{jt}+\epsilon_{jt}^2\right)-(\gamma_i^2\sigma_1^2+\sigma^2)(\gamma_j^2\sigma_1^2+\sigma^2)\right)\\
	&=& O(\frac{1}{nT})+\frac{1}{n^2 T^2}\sum_{j\neq i}^{n}\sum_{t=1}^T\sum_{i=1}^{n}
	\left(\gamma_i^2\gamma_j^2 Ef_t^4+\gamma_i^2 Ef_t^2 E\epsilon_{jt}^2+\gamma_j^2 Ef_t^2 E\epsilon_{it}^2+E\epsilon_{it}^2E\epsilon_{jt}^2-(\gamma_i^2\sigma_1^2+\sigma^2)(\gamma_j^2\sigma_1^2+\sigma^2)\right)\\
	&=& O(\frac{1}{nT})+\frac{1}{n^2 T^2}\sum_{j\neq i}^{n}\sum_{t=1}^T\sum_{i=1}^{n}
	\left(\gamma_i^2\gamma_j^2 Ef_t^4+\gamma_i^2 \sigma_1^2\sigma^2+\gamma_j^2 \sigma_1^2\sigma^2+\sigma^4-(\gamma_i^2\sigma_1^2+\sigma^2)(\gamma_j^2\sigma_1^2+\sigma^2)\right)\\
	&=& O(\frac{1}{nT})+\frac{1}{n^2 T^2}\sum_{j\neq i}^{n}\sum_{t=1}^T\sum_{i=1}^{n}
	\left(\gamma_i^2\gamma_j^2 Ef_t^4-\gamma_i^2\gamma_j^2\sigma_1^4\right)\\
	&\leq& O(\frac{1}{nT})+\frac{1}{n^2 T}\sum_{j=1}^{n}\gamma_j^2\sum_{i=1}^{n}\gamma_i^2
	\left(Ef_1^4-\sigma_1^4\right)\\
	&=& O(\frac{1}{nT})+O(\frac{1}{n^2 T}).
\end{eqnarray*}

Proof of part (b).
\begin{eqnarray*}
	\frac{1}{n}tr(S^2)
	&=&\frac{1}{n}tr\left[\left(\frac{1}{T}\sum_{t=1}^T v_t v_t^\prime\right)\left(\frac{1}{T}\sum_{s=1}^T v_s v_s^\prime\right)\right]\\
	&=&\frac{1}{nT^2}\sum_{t=1}^T\sum_{i=1}^{n} v_{it}^4
	+\frac{1}{nT^2}\sum_{t=1}^T\sum_{s\neq t}^T\sum_{i=1}^{n} v_{it}^2 v_{is}^2\\
	&&+\frac{1}{nT^2}\sum_{t=1}^T\sum_{i=1}^{n}\sum_{j\neq i}^{n} v_{it}^2 v_{jt}^2
	+\frac{1}{nT^2}\sum_{t=1}^T\sum_{s\neq t}^T\sum_{i=1}^{n}\sum_{j\neq i}^{n} v_{it} v_{is} v_{jt} v_{js}\\
	&=&O_p\left(\frac{1}{T}\right)+\left[\sigma^4+O_p\left(\frac{1}{\sqrt{T}}\right)\right]+\left[\frac{n}{T}\sigma^4+O_p\left(\frac{1}{\sqrt{T}}\right)\right]+O_p\left(\frac{1}{T}\right).
\end{eqnarray*}

This uses the following four results.

(1) $\frac{1}{nT^2}\sum_{t=1}^T\sum_{i=1}^{n} v_{it}^4=O_p\left(\frac{1}{T}\right)$ due to
\begin{eqnarray*}
	E\left(\frac{1}{nT^2}\sum_{t=1}^T\sum_{i=1}^{n} v_{it}^4\right)^2
	&=& \frac{1}{n^2 T^4}\sum_{t=1}^T\sum_{i=1}^{n}E v_{it}^8
	+\frac{1}{n^2 T^4}\sum_{j\neq i}^{n}\sum_{t=1}^T\sum_{i=1}^{n}E v_{it}^4 v_{jt}^4\\
	&&+\frac{1}{n^2 T^4}\sum_{s\neq t}^{n}\sum_{t=1}^T\sum_{i=1}^{n}E v_{it}^4 v_{is}^4
	+\frac{1}{n^2 T^4}\sum_{s\neq t}^T\sum_{j\neq i}^{n}\sum_{t=1}^T\sum_{i=1}^{n}E v_{it}^4 v_{js}^4\\
	&=& O\left(\frac{1}{nT^3}\right)+O\left(\frac{1}{T^3}\right)+O\left(\frac{1}{T^3}\right)+O\left(\frac{1}{T^2}\right).
\end{eqnarray*}

(2) 
\begin{eqnarray*}
	&&\frac{1}{nT^2}\sum_{t=1}^T\sum_{s\neq t}^T\sum_{i=1}^{n} v_{it}^2 v_{is}^2\\
	&=&\frac{1}{nT^2}\sum_{t=1}^T\sum_{s\neq t}^T\sum_{i=1}^{n} Ev_{it}^2 Ev_{is}^2+\frac{1}{nT^2}\sum_{t=1}^T\sum_{s\neq t}^T\sum_{i=1}^{n} \left(v_{it}^2 v_{is}^2-Ev_{it}^2 Ev_{is}^2\right)\\
	&=&\frac{1}{nT^2}\sum_{t=1}^T\sum_{s\neq t}^T\sum_{i=1}^{n}(\sigma^2+\sigma^2he_i^2)^2+\frac{1}{nT^2}\sum_{t=1}^T\sum_{s\neq t}^T\sum_{i=1}^{n} \left(v_{it}^2 v_{is}^2-Ev_{it}^2 Ev_{is}^2\right)\\
	&=&\left(\frac{T-1}{T}\sigma^4+\frac{T-1}{nT}2\sigma^4h+\frac{T-1}{nT}\sigma^4h^2\sum_{i=1}^n e_i^4\right) +\frac{1}{nT^2}\sum_{t=1}^T\sum_{s\neq t}^T\sum_{i=1}^{n} \left(v_{it}^2 v_{is}^2-Ev_{it}^2 Ev_{is}^2\right)\\
	&=&\left[\frac{T-1}{T}\sigma^4+O_p\left(\frac{1}{n}\right)\right]+\left[O_p\left(\frac{1}{\sqrt{T}}\right)\right]\\
	&=&\sigma^4+O_p\left(\frac{1}{\sqrt{T}}\right)+O_p\left(\frac{1}{n}\right),
\end{eqnarray*}
due to
\begin{eqnarray*}
	&&E\left(\frac{1}{nT^2}\sum_{t=1}^T\sum_{s\neq t}^T\sum_{i=1}^{n} \left(v_{it}^2 v_{is}^2-Ev_{it}^2 Ev_{is}^2\right)\right)^2\\
	&=&\frac{1}{n^2T^4}\sum_{t_1=1}^T\sum_{s_1\neq t_1}^T\sum_{i_1=1}^{n}\sum_{t_2=1}^T\sum_{s_2\neq t_2}^T\sum_{i_2=1}^{n}E\left(v_{i_1t_1}^2 v_{i_1s_1}^2-Ev_{i_1t_1}^2 Ev_{i_1s_1}^2\right)\left(v_{i_2t_2}^2 v_{i_2s_2}^2-Ev_{i_2t_2}^2 Ev_{i_2s_2}^2\right)\\
	&=&\frac{1}{n^2T^4}\sum_{t_1=1}^T\sum_{s_1\neq t_1}^T\sum_{i_1=1}^{n}\sum_{t_2=1}^T\sum_{s_2\neq t_2}^T\sum_{i_2=1}^{n}\cov(v_{i_1t_1}^2 v_{i_1s_1}^2,v_{i_2t_2}^2 v_{i_2s_2}^2)\\
	&=&\frac{1}{n^2T^4}\sum_{t_1=1}^T\sum_{s_1\neq t_1}^T\sum_{i_1=1}^{n}\sum_{t_2\neq t_1,s_1}^T\sum_{s_2\neq t_2,t_1,s_1}^T\sum_{i_2=1}^{n}\cov(v_{i_1t_1}^2 v_{i_1s_1}^2,v_{i_2t_2}^2 v_{i_2s_2}^2)\\ &&+\frac{1}{n^2T^4}\sum_{t_1=1}^T\sum_{s_1\neq t_1}^T\sum_{i_1=1}^{n}\sum_{s_2\neq t_1,s_1}^T\sum_{i_2=1}^{n}\cov(v_{i_1t_1}^2 v_{i_1s_1}^2,v_{i_2t_1}^2 v_{i_2s_2}^2)\\
	&&+\frac{1}{n^2T^4}\sum_{t_1=1}^T\sum_{s_1\neq t_1}^T\sum_{i_1=1}^{n}\sum_{s_2\neq t_1,s_1}^T\sum_{i_2=1}^{n}\cov(v_{i_1t_1}^2 v_{i_1s_1}^2,v_{i_2s_1}^2 v_{i_2s_2}^2)\\
	&&+\frac{1}{n^2T^4}\sum_{t_1=1}^T\sum_{s_1\neq t_1}^T\sum_{i_1=1}^{n}\sum_{t_2\neq t_1,s_1}^T\sum_{i_2=1}^{n}\cov(v_{i_1t_1}^2 v_{i_1s_1}^2,v_{i_2t_2}^2 v_{i_2t_1}^2)\\
	&&+\frac{1}{n^2T^4}\sum_{t_1=1}^T\sum_{s_1\neq t_1}^T\sum_{i_1=1}^{n}\sum_{t_2\neq t_1,s_1}^T\sum_{i_2=1}^{n}\cov(v_{i_1t_1}^2 v_{i_1s_1}^2,v_{i_2t_2}^2 v_{i_2s_1}^2)\\
	&=& 0+O(\frac{1}{T})+O(\frac{1}{T})+O(\frac{1}{T})+O(\frac{1}{T})
	= O(\frac{1}{T}).
\end{eqnarray*}

Proof of part (c) follows directly from \cite{baltagi2017asymptotic}.

Proof of part (d) follows direclt from Lemma \ref{lemma_a_alternative} below.
\begin{lemma}\label{lemma_a_alternative}
	Under Assumptions 1, 2 and the weak factor alternative, 
	
	(a) $\frac{1}{n}tr(SA_1)=O_p(\frac{1}{T^2})+O_p(\frac{1}{nT})+O_p(\frac{1}{T\sqrt{nT}})$;
	
	(b) $\frac{1}{n}tr(SA_3)=O_p(\frac{1}{nT})$;
	
	(c) $\frac{1}{n}tr(A_1^2)=O_p(\frac{1}{nT^2})$;
	
	(d) $\frac{1}{n}tr(A_1A_2)=O_p(\frac{1}{T^2})$;
	
	(e) $\frac{1}{n}tr(A_1A_3)=O_p(\frac{1}{nT^2})$;
	
	(f) $\frac{1}{n}tr(A_3^2)=O_p(\frac{1}{nT^2})$;
	
	(g) $\frac{1}{n}tr(SA_0)=\frac{1}{T}\sigma^4+\frac{n}{T^2}\sigma^4+O_p(\frac{1}{T\sqrt{T}})$; 
	
	(h) $\frac{1}{n}tr(A_0^2)=\frac{n}{T^2}\sigma^4+O_p(\frac{\sqrt{n}}{T^2})$;
	
	(i) $\frac{1}{n}tr(A_0 A_1)=O_p(\frac{1}{T^2})$;
	
	(j) $\frac{1}{n}tr(A_0 A_3)=O_p(\frac{1}{nT^2})$.
\end{lemma}

This lemma can be proved following the same proof line of Lemma 3 in \cite{baltagi2017asymptotic}.

Then, the first result in Proposition \ref{diff_between_U_weak}, the difference between John's test statistic $U$ and the residual-based John's test statistic $\hat{U}$, is calculated. Define $W_1=\frac{1}{n}tr(\hat{S})-\frac{1}{n}tr(S)$, $W_2=\frac{1}{n}tr(\hat{S}^2)-\frac{1}{n}tr(S^2)$
\begin{eqnarray*}
	T(\hat{U}-U)&=&\frac{TW_2(\frac{1}{n}trS)^2-2TW_1\frac{1}{n}trS\frac{1}{n}tr S^2-T W_1^2\frac{1}{n}tr S^2}{\left(\frac{1}{n}trS+W_1\right)^2\left(\frac{1}{n}trS\right)^2}\\
	&=&\frac{c_T\sigma^8+O_p(\frac{1}{\sqrt{T}})}{\sigma^8+O_p(\frac{1}{\sqrt{nT}})+O_p(\frac{1}{T})}\\
	&=& c_T+O_p\left(\frac{1}{\sqrt{T}}\right).
\end{eqnarray*}

To prove the second result in Proposition \ref{diff_between_U_weak}, it is enough to prove that $E|\hat{\nu}_{it}|^4 = E(\nu_{it}^4) +O_p(1/T)$. For notation simplicity, the proof given here is for $r=1$. Using $\left(\sum_{j=1}^r h_j \right)^2 \leq r \sum_{j=1}^r h_j^2$ repeatedly, the case where $r>1$ can be proved similarly, as long as $r$ is fixed. As $\nu_{it} = \lambda_{i1}f_{t1}+\epsilon_{it}$, we have $E(\nu_{it}^2 )=h/n+1$. As $h\to d \in (0,\infty)$ and $n\to \infty$, we have
\begin{eqnarray*}
	E|\hat{\nu}_{it}|^4 &=& E\left|\nu_{it} - \bar{\nu}_{i.}-\tilde{x}_{it}^\prime (\hat{\beta}-\beta)\right|^4  \\
	&=& E(\nu_{it}-\bar{\nu}_{i.})^4 + E\left[\tilde{x}_{it}^\prime (\hat{\beta}-\beta)\right]^4 + 6 E (\nu_{it}-\bar{\nu}_{i.})^2\left[\tilde{x}_{it}^\prime (\hat{\beta}-\beta)\right]^2,
\end{eqnarray*}
where
\begin{eqnarray*}
	E(\nu_{it}-\bar{\nu}_{i.})^4 &=& E \left(\nu_{it}^4 + \bar{\nu}_{i.}^4 -4 \nu_{it}^3\bar{\nu}_{i.} + 6 \nu_{it}^2\bar{\nu}_{i.}^2 -4\nu_{it}\bar{\nu}_{i.}^3  \right) \\
	&=& E(\nu_{it}^4)\left(1-\frac{4}{T}+\frac{6}{T^2}-\frac{3}{T^3} \right) + \left(\frac{h}{n}+1 \right)^2\left(\frac{6(T-1)}{T^2} -\frac{11(T-1)}{T^3} \right)\\
	&=& E(\nu_{it}^4) +O_p\left(\frac{1}{T}\right),\\
	&&E\left[\tilde{x}_{it}^\prime (\hat{\beta}-\beta)\right]^4 = O_p\left(\frac{1}{n^2T^2} \right),\\
	&&E (\nu_{it}-\bar{\nu}_{i.})^2\left[\tilde{x}_{it}^\prime (\hat{\beta}-\beta)\right]^2 = O_p\left(\frac{1}{nT} \right).
\end{eqnarray*}
Therefore,
\begin{eqnarray*}
	E|\hat{\nu}_{it}|^4 = E(\nu_{it}^4) +O_p\left(\frac{1}{T} \right).
\end{eqnarray*}
The proof of Proposition \ref{diff_between_U_weak} is complete. 

\section{Proof of Theorem \ref{ultra_high_weak}}

The result of Proposition \ref{U_ulpa} is obtained directly from Theorem 3.2 of \cite{li2016testing}. The proof of Proposition \ref{conjection_general_ulpa} is the same as that of Proposition \ref{prop3}.

\section{Proof of Theorem \ref{J_H1}}
The proof of Theorem \ref{J_H1} is based on Theorem 2.2 of \cite{yin2021spectral}, from which we have the joint distribution of $n^{-1} tr(S_T)$ and $n^{-1} tr(S_T^2)$,
\[T\left(\begin{array}{c}
	n^{-1} tr(S_T)-\frac{m_1}{n}\\
	n^{-1} tr(S_T^2)-\frac{m_2}{n}
\end{array} \right) \stackrel{D}{\rightarrow} N\left( \left(\begin{array}{c}
	0 \\ 0
\end{array} \right), \frac{1}{c^2}\cdot \boldsymbol{\Psi}_2 \right), \]
where $$
m_1=\operatorname{tr}(\Sigma_n), \quad
 m_2=T^{-1}\left[(\gamma_4-3) \operatorname{tr}(\Sigma_n \circ \Sigma_n)+\operatorname{tr}^{2}(\Sigma_n)+(T+1) \operatorname{tr}\left(\Sigma_n^{2}\right)\right],
$$
and the covariance matrix is $\boldsymbol{\Psi}_{2}=\left(\psi_{i j}\right)_{2 \times 2}$ with its entries
$$
\begin{aligned}
	\psi_{11}=& T^{-1}\left[(\gamma_4-3) \operatorname{tr}\left(\Sigma_n \circ \Sigma_n \right)+2 \operatorname{tr}\left(\Sigma_n^{2}\right)\right], \\
	\psi_{12}=& \psi_{21}=T^{-2}\left[4 \operatorname{tr}\left(\Sigma_n^{2}\right) \operatorname{tr}(\Sigma_n)+2 (\gamma_4-3) \operatorname{tr}\left(\Sigma_n \circ \Sigma_n \right) \operatorname{tr} \Sigma_n+2 (\gamma_4-3) T \operatorname{tr}\left(\Sigma_n \circ \Sigma_n^{2}\right)+4 T \operatorname{tr}\left(\Sigma_n^{3}\right)\right], \\
	\psi_{22}=& T^{-3}\left[8 \operatorname{tr}\left(\Sigma_n^{2}\right) \operatorname{tr}^{2}(\Sigma_n)+4 (\gamma_4-3) \operatorname{tr}^{2}(\Sigma_n) \operatorname{tr}\left( \Sigma_n \circ \Sigma_n\right)+16 T \operatorname{tr}(\Sigma_n) \operatorname{tr}\left(\Sigma_n^{3}\right)\right.\\
	&\left.+4 T \operatorname{tr}^{2}\left(\Sigma_n^{2}\right)+8 (\gamma_4-3)T \operatorname{tr}\left(\Sigma_n \circ \Sigma_n^{2}\right) \operatorname{tr} \Sigma_n+4 (\gamma_4-3) T^{2} \operatorname{tr}\left(\Sigma_n^{2} \circ \Sigma_n^{2}\right)+8 T^{2} \operatorname{tr}\left(\Sigma_n^{4}\right)\right].
\end{aligned}
$$

Define the function $q(x, y)=\frac{y}{x^2}-1$, then $U=q\(n^{-1} tr(S_T), n^{-1} tr(S_T^2)\)$.

We have
\begin{eqnarray*}
	&&\frac{\partial q}{\partial x}\(\frac{m_1}{n},\frac{m_2}{n} \)=-\frac{2n^2 m_2}{m_1^3},\\
	&&\frac{\partial q}{\partial y}\(\frac{m_1}{n},\frac{m_2}{n} \)= \frac{n^2}{m_1^2},\\
	&& q \(\frac{m_1}{n},\frac{m_2}{n} \) = \frac{n m_2}{m_1^2}-1.
\end{eqnarray*}
By the delta method,
\[n\(U-q\(\frac{m_1}{n},\frac{m_2}{n}\) \)\stackrel{D}{\rightarrow} N(0, \sigma^2),  \]
where
\begin{eqnarray*}
	\sigma^2 \sim \(\begin{array}{c}
		\frac{\partial q}{\partial x}\(\frac{m_1}{n},\frac{m_2}{n} \) \\ \frac{\partial q}{\partial y}\(\frac{m_1}{n},\frac{m_2}{n} \)
	\end{array} \)^T \cdot \(\frac{1}{c^2} \boldsymbol{\Psi}_2 \)\(\begin{array}{c}
		\frac{\partial q}{\partial x}\(\frac{m_1}{n},\frac{m_2}{n} \) \\ \frac{\partial q}{\partial y}\(\frac{m_1}{n},\frac{m_2}{n} \)
	\end{array} \).
\end{eqnarray*}
Therefore,
\[T\( U- \frac{T m_2}{(m_1)^2}+1\) \stackrel{D}{\rightarrow} N(0,\sigma^2). \]
The proof of Theorem \ref{J_H1} is complete.

\section{Proof of Proposition \ref{diff_between_U_diverge_factor}}
Following the proof of Proposition 1 in \cite{baltagi2017asymptotic}, the proof is conducted as follows. But we don't calculate $\frac{1}{n}tr(SA_0)$ and $W_1=\frac{1}{n}tr(\hat{S})-\frac{1}{n}tr(S)$ explictly. Instead, we represent them as a function of $\frac{1}{n}trS$. We could see that those two items are canceled in the following proof.

Define
\begin{eqnarray*}
&&R_1=\frac{1}{nT^3}\left(2\sum_{t=1}^{T}\sum_{s\neq t}^{T}\sum_{i=1}^{n}\sum_{j=1}^{n}v_{jt}^2v_{it}v_{is}+\sum_{t=1}^{T}\sum_{s\neq t}^{T}\sum_{\tau\neq s,t}^{T}\sum_{i=1}^{n}\sum_{j=1}^{n}v_{it}v_{is}v_{j\tau}v_{jt}\right),\\
&&R_2=-\frac{1}{nT^2}\sum_{t=1}^{T}\sum_{s\neq t}^{T}\sum_{i=1}^{n}v_{is}v_{it}.
\end{eqnarray*}

\begin{lemma}
Under Assumptions 1,2 and the divergent number of factors alternative,

(a) $\frac{1}{n}tr(SA_1)=O_p\left(\frac{r}{nT\sqrt{T}}\right)+O_p\left(\frac{\sqrt{n}}{T^2}\right)+O_p\left(\frac{1}{T\sqrt{n}}\right)$;

(b) $\frac{1}{n}tr(SA_3)=O_p\left(\frac{r}{nT}\right)$;

(c) $\frac{1}{n}tr(A_1^2)=O_p\left(\frac{1}{T^2}\right)$;

(d) $\frac{1}{n}tr(A_1A_2)=O_p\left(\frac{1}{T^2}\right)$;

(e) $\frac{1}{n}tr(A_1A_3)=O_p\left(\frac{1}{\sqrt{n}T^2}\right)$;

(f) $\frac{1}{n}tr(A_3^2)=O_p\left(\frac{1}{nT^2}\right)$;

(g) $\frac{1}{n}tr(SA_0)=\frac{1}{nT}trS^2+R_1$;

(h) $\frac{1}{n}tr(A_0^2)=\frac{n}{T^2}\left(\frac{1}{n}trS\right)^2+O_p\left(\frac{\sqrt{n}}{T^2}\right)$;

(i) $\frac{1}{n}tr(A_0 A_1)=O_p\left(\frac{\sqrt{n}}{T^2}\right)$;

(j) $\frac{1}{n}tr(A_0 A_3)=O_p\left(\frac{1}{T^2}\right)$.

\end{lemma}
The proofs of parts $(c), (d), (e), (f), (i), (j)$ are trivial, and we only prove parts $(a),(b),(g),(h)$.

\subsubsection*{Proof of part (a)}
\begin{eqnarray}\label{f.1}
	\frac{1}{n}tr(SA_1) &=& \frac{1}{n}tr\left(S\frac{1}{T}\sum_{t=1}^{T}\tilde{x}_t\left(\hat{\beta}-\beta\right)\tilde{v}_t^{\prime}\right)\nonumber\\
	&=& \frac{1}{n}tr\left(\frac{1}{T}\sum_{s=1}^{T}v_s v_s^{\prime}\frac{1}{T}\sum_{t=1}^{T}\tilde{x}_t\left(\hat{\beta}-\beta\right)\tilde{v}_t^{\prime}\right)\nonumber\\
		&=& \frac{1}{n T^2}tr\left(\sum_{s=1}^{T}\sum_{t=1}^{T}\tilde{v}_t^{\prime}v_s v_s^{\prime}\tilde{x}_t\left(\hat{\beta}-\beta\right)\right)\nonumber\\
	&=& \frac{1}{n T^2}\sum_{s=1}^{T}\sum_{t=1}^{T}\tilde{v}_t^{\prime}v_s v_s^{\prime}\tilde{x}_t\left(\hat{\beta}-\beta\right)\nonumber\\
	&=& \left(\frac{1}{n T^2}\sum_{s=1}^{T}\sum_{t=1}^{T}v_t^{\prime}v_s v_s^{\prime}\tilde{x}_t-\frac{1}{n T^2}\sum_{s=1}^{T}\sum_{t=1}^{T}\bar{v}^{\prime}v_s v_s^{\prime}\tilde{x}_t\right)\left(\hat{\beta}-\beta\right).
\end{eqnarray}
Consider the first term of (\ref{f.1}),
\begin{eqnarray*}
	\frac{1}{n T^2}\sum_{s=1}^{T}\sum_{t=1}^{T}v_t^{\prime}v_s v_s^{\prime}\tilde{x}_t
	&=&\frac{1}{n T^2}\sum_{t=1}^{T}\sum_{i=1}^{n}v_{it}^3\tilde{x}_{it}^{\prime}+\frac{1}{n T^2}\sum_{t=1}^{T}\sum_{j\neq i}^{n}\sum_{i=1}^{n}v_{it}^2 v_{jt}\tilde{x}_{jt}^{\prime}\\
	&& +\frac{1}{n T^2}\sum_{s\neq t}^{T}\sum_{t=1}^{T}\sum_{i=1}^{n}v_{it}v_{is}^2\tilde{x}_{it}^{\prime}+\frac{1}{n T^2}\sum_{s\neq t}^{T}\sum_{t=1}^{T}\sum_{j\neq i}^{n}\sum_{i=1}^{n}v_{it}v_{is}v_{js}\tilde{x}_{jt}^{\prime}.
\end{eqnarray*}
We examine the orders of the last three terms, respectively. First,
\begin{eqnarray*}
	&&E\left(\frac{1}{n T^2}\sum_{t=1}^{T}\sum_{j\neq i}^{n}\sum_{i=1}^{n}v_{it}^2 v_{jt}\tilde{x}_{jt}\right)^2\\
	&=&\frac{1}{n^2T^4}\sum_{t_1=1}^{T}\sum_{j_1\neq i_1}^{n}\sum_{i_1=1}^{n}\sum_{t_2=1}^{T}\sum_{j_2\neq i_2}^{n}\sum_{i_2=1}^{n}Ev_{i_1t_1}^2 v_{j_1t_1}\tilde{x}_{j_1t_1}v_{i_2t_2}^2 v_{j_2t_2}\tilde{x}_{j_2t_2}\\
	&=&\frac{1}{n^2T^4}\sum_{t_1=1}^{T}\sum_{j_1\neq i_1}^{n}\sum_{i_1=1}^{n}\sum_{j_2\neq i_2}^{n}\sum_{i_2=1}^{n}Ev_{i_1t_1}^2 v_{j_1t_1}v_{i_2t_1}^2 v_{j_2t_1}E\tilde{x}_{j_1t_1}\tilde{x}_{j_2t_1}\\
	&&+\frac{1}{n^2T^4}\sum_{t_1=1}^{T}\sum_{j_1\neq i_1}^{n}\sum_{i_1=1}^{n}\sum_{t_2\neq t_1}^{T}\sum_{j_2\neq i_2}^{n}\sum_{i_2=1}^{n}Ev_{i_1t_1}^2 v_{j_1t_1}Ev_{i_2t_2}^2 v_{j_2t_2}E\tilde{x}_{j_1t_1}\tilde{x}_{j_2t_2}\\
	&=&O\left(\frac{n^2}{T^3}\right)
	+\frac{1}{n^2T^4}\sum_{t_1=1}^{T}\sum_{j_1\neq i_1}^{n}\sum_{i_1=1}^{n}\sum_{t_2\neq t_1}^{T}\sum_{j_2\neq i_2}^{n}\sum_{i_2=1}^{n}Ev_{i_1t_1}^2 v_{j_1t_1}Ev_{i_2t_2}^2 v_{j_2t_2}E\tilde{x}_{j_1t_1}\tilde{x}_{j_2t_2}\\
	&\leq&O\left(\frac{n^2}{T^3}\right)
	+\frac{\sqrt{C_1}}{n^2T^4}\sum_{t_1=1}^{T}\sum_{j_1\neq i_1}^{n}\sum_{i_1=1}^{n}\sum_{t_2\neq t_1}^{T}\sum_{j_2\neq i_2}^{n}\sum_{i_2=1}^{n}Ev_{i_1t_1}^2 v_{j_1t_1}Ev_{i_2t_2}^2 v_{j_2t_2},
\end{eqnarray*}
by plugging in
\begin{eqnarray*}
	Ev_{i_1t_1}^2v_{j_1t_1}&=&E\left(\sum_{k_1=1}^{r}\gamma_{i_1k_1}f_{t_1k_1}+\epsilon_{i_1t_1}\right)^2\left(\sum_{k_2=1}^{r}\gamma_{j_1k_2}f_{t_1k_2}+\epsilon_{j_1t_1}\right)\\
	&=&E\left(\sum_{k_1=1}^{r}\gamma_{i_1k_1}f_{t_1k_1}+\epsilon_{i_1t_1}\right)^2\left(\sum_{k_2=1}^{r}\gamma_{j_1k_2}f_{t_1k_2}\right)\\
	&=&E\left(\left(\sum_{k_1=1}^{r}\gamma_{i_1k_1}f_{t_1k_1}\right)^2+2\epsilon_{i_1t_1}\left(\sum_{k_1=1}^{r}\gamma_{i_1k_1}f_{t_1k_1}\right)+\epsilon_{i_1t_1}^2\right)\left(\sum_{k_2=1}^{r}\gamma_{j_1k_2}f_{t_1k_2}\right)\\
	&=&E\left(\sum_{k_1=1}^{r}\gamma_{i_1k_1}f_{t_1k_1}\right)^2\left(\sum_{k_2=1}^{r}\gamma_{j_1k_2}f_{t_1k_2}\right)\\
	&=&\sum_{k=1}^{r}\gamma_{i_1k}^2\gamma_{j_1k}Ef_{t_1k}^3,
\end{eqnarray*}
and
\begin{eqnarray*}
	&&\frac{1}{n^2T^4}\sum_{t_1=1}^{T}\sum_{j_1\neq i_1}^{n}\sum_{i_1=1}^{n}\sum_{t_2\neq t_1}^{T}\sum_{j_2\neq i_2}^{n}\sum_{i_2=1}^{n}Ev_{i_1t_1}^2 v_{j_1t_1}Ev_{i_2t_2}^2 v_{j_2t_2}\\
	&=&\frac{1}{n^2T^4}\sum_{t_1=1}^{T}\sum_{j_1\neq i_1}^{n}\sum_{i_1=1}^{n}\sum_{t_2\neq t_1}^{T}\sum_{j_2\neq i_2}^{n}\sum_{i_2=1}^{n}\left(\sum_{k_1=1}^{r}\gamma_{i_1k_1}^2\gamma_{j_1k_1}Ef_{t_1k_1}^3\right)\left(\sum_{k_2=1}^{r}\gamma_{i_2k_2}^2\gamma_{j_2k_2}Ef_{t_1k_2}^3\right)\\
	&=&\frac{1}{n^2T^4}\sum_{t_1=1}^{T}\sum_{t_2\neq t_1}^{T}\left(\sum_{k_1=1}^{r}\sum_{j_1\neq i_1}^{n}\sum_{i_1=1}^{n}\gamma_{i_1k_1}^2\gamma_{j_1k_1}Ef_{t_1k_1}^3\right)\left(\sum_{k_2=1}^{r}\sum_{j_2\neq i_2}^{n}\sum_{i_2=1}^{n}\gamma_{i_2k_2}^2\gamma_{j_2k_2}Ef_{t_1k_2}^3\right)\\
	&=&\frac{1}{n^2T^4}O(T^2nr^2)=O\left(\frac{r^2}{nT^2}\right),
\end{eqnarray*}
we have, $E\left(\frac{1}{n T^2}\sum_{t=1}^{T}\sum_{j\neq i}^{n}\sum_{i=1}^{n}v_{it}^2 v_{jt}\tilde{x}_{jt}\right)^2=O\left(\frac{n^2}{T^3}\right)+O\left(\frac{r^2}{nT^2}\right)$.

Second,
\begin{eqnarray*}
	&&E\left(\frac{1}{n T^2}\sum_{s\neq t}^{T}\sum_{t=1}^{T}\sum_{i=1}^{n}v_{it}v_{is}^2\tilde{x}_{it}^{\prime}\right)^2\\
	&=&\frac{1}{n^2 T^4}\sum_{s_1\neq t_1}^{T}\sum_{t_1=1}^{T}\sum_{i_1=1}^{n}\sum_{s_2\neq t_2}^{T}\sum_{t_2=1}^{T}\sum_{i_2=1}^{n}Ev_{i_1t_1}v_{i_1s_1}^2v_{i_2t_2}v_{i_2s_2}^2E\tilde{x}_{i_1t_1}\tilde{x}_{i_2t_2}\\
	&=&\frac{1}{n^2 T^4}\sum_{s_1\neq t_1}^{T}\sum_{t_1=1}^{T}\sum_{i_1=1}^{n}\sum_{s_2\neq t_2}^{T}\sum_{i_2=1}^{n}Ev_{i_1t_1}v_{i_1s_1}^2v_{i_2t_1}v_{i_2s_2}^2E\tilde{x}_{i_1t_1}\tilde{x}_{i_2t_1}\\
	&&+\frac{1}{n^2 T^4}\sum_{s_1\neq t_1}^{T}\sum_{t_1=1}^{T}\sum_{i_1=1}^{n}\sum_{s_2\neq t_2}^{T}\sum_{t_2\neq t_1}^{T}\sum_{i_2=1}^{n}Ev_{i_1t_1}v_{i_1s_1}^2v_{i_2t_2}v_{i_2s_2}^2E\tilde{x}_{i_1t_1}\tilde{x}_{i_2t_2}\\
	&=&\frac{1}{n^2 T^4}\sum_{s_1\neq t_1}^{T}\sum_{t_1=1}^{T}\sum_{i_1=1}^{n}\sum_{s_2\neq t_2}^{T}\sum_{i_2=1}^{n}Ev_{i_1t_1}v_{i_1s_1}^2v_{i_2t_1}v_{i_2s_2}^2E\tilde{x}_{i_1t_1}\tilde{x}_{i_2t_1}+0\\
	&=&O\left(\frac{1}{T}\right).
\end{eqnarray*}
Third,
\begin{eqnarray*}
	&&E\left(\frac{1}{n T^2}\sum_{s\neq t}^{T}\sum_{t=1}^{T}\sum_{j\neq i}^{n}\sum_{i=1}^{n}v_{it}v_{is}v_{js}\tilde{x}_{jt}\right)^2\\
	&=&\frac{1}{n^2 T^4}\sum_{s_1\neq t}^{T}\sum_{s_2\neq t}^{T}\sum_{t=1}^{T}\sum_{j_1\neq i_1}^{n}\sum_{i_1=1}^{n}\sum_{j_2\neq i_2}^{n}\sum_{i_2=1}^{n}Ev_{i_1t}v_{i_2t}Ev_{i_1s_1}v_{j_1s_1}v_{i_2s_2}v_{j_2s_2}E\tilde{x}_{j_1t}\tilde{x}_{j_2t}\\
	&\leq&\frac{\sqrt{C_1}}{n^2 T^4}\sum_{s_1\neq t}^{T}\sum_{s_2\neq t}^{T}\sum_{t=1}^{T}\sum_{j_1\neq i_1}^{n}\sum_{i_1=1}^{n}\sum_{j_2\neq i_2}^{n}\sum_{i_2=1}^{n}Ev_{i_1t}v_{i_2t}Ev_{i_1s_1}v_{j_1s_1}v_{i_2s_2}v_{j_2s_2}\\
	&=&\frac{\sqrt{C_1}}{n^2 T^4}\sum_{s_1\neq t}^{T}\sum_{t=1}^{T}\sum_{j_1\neq i_1}^{n}\sum_{i_1=1}^{n}\sum_{j_2\neq i_2}^{n}\sum_{i_2=1}^{n}Ev_{i_1t}v_{i_2t}Ev_{i_1s_1}v_{j_1s_1}v_{i_2s_1}v_{j_2s_1}\\
	&&+\frac{\sqrt{C_1}}{n^2 T^4}\sum_{s_1\neq t}^{T}\sum_{s_2\neq s_1,t}^{T}\sum_{t=1}^{T}\sum_{j_1\neq i_1}^{n}\sum_{i_1=1}^{n}\sum_{j_2\neq i_2}^{n}\sum_{i_2=1}^{n}Ev_{i_1t}v_{i_2t}Ev_{i_1s_1}v_{j_1s_1}Ev_{i_2s_2}v_{j_2s_2}\\
	&=&O\left(\frac{r}{nT}\right)+\frac{\sqrt{C_1}}{n^2 T^4}\sum_{s_1\neq t}^{T}\sum_{s_2\neq s_1,t}^{T}\sum_{t=1}^{T}\sum_{j_1\neq i_1}^{n}\sum_{i_1=1}^{n}\sum_{j_2\neq i_2}^{n}\sum_{i_2=1}^{n}\Sigma_{i_1i_2}\Sigma_{i_1j_1}\Sigma_{i_2j_2}\\
	&=&O\left(\frac{r}{nT}\right)+\frac{\sqrt{C_1}}{n^2 T^4}\sum_{s_1\neq t}^{T}\sum_{s_2\neq s_1,t}^{T}\sum_{t=1}^{T}\sum_{j_1\neq i_1}^{n}\sum_{i_1=1}^{n}\sum_{j_2\neq i_2}^{n}\Sigma_{i_1i_1}\Sigma_{i_1j_1}\Sigma_{i_2j_2}\\
	&&+\frac{\sqrt{C_1}}{n^2 T^4}\sum_{s_1\neq t}^{T}\sum_{s_2\neq s_1,t}^{T}\sum_{t=1}^{T}\sum_{j_1\neq i_1}^{n}\sum_{i_1=1}^{n}\sum_{j_2\neq i_2}^{n}\sum_{i_2\neq i_1}^{n}\Sigma_{i_1i_2}\Sigma_{i_1j_1}\Sigma_{i_2j_2}\\
	&=&O\left(\frac{r}{nT}\right)\\
	&&+\frac{\sqrt{C_1}}{n^2 T^4}\sum_{s_1\neq t}^{T}\sum_{s_2\neq s_1,t}^{T}\sum_{t=1}^{T}\sum_{j_1\neq i_1}^{n}\sum_{i_1=1}^{n}\sum_{j_2\neq i_1}^{n}\left(1+\sum_{k_1=1}^{r}h_{k_1}e_{i_1k_1}^2\right)\left(\sum_{k_2=1}^{r}h_{k_2}e_{i_1k_2}e_{j_1k_2}\right)\left(\sum_{k_3=1}^{r}h_{k_3}e_{i_1k_3}e_{j_2k_3}\right)\\
	&&+\frac{\sqrt{C_1}}{n^2 T^4}\sum_{s_1\neq t}^{T}\sum_{s_2\neq s_1,t}^{T}\sum_{t=1}^{T}\sum_{j_1\neq i_1}^{n}\sum_{i_1=1}^{n}\sum_{j_2\neq i_2}^{n}\sum_{i_2\neq i_1}^{n}\left(\sum_{k_1=1}^{r}h_{k_1}e_{i_1k_1}e_{i_2k_1}\right)\left(\sum_{k_2=1}^{r}h_{k_2}e_{i_1k_2}e_{j_1k_2}\right)\left(\sum_{k_3=1}^{r}h_{k_3}e_{i_2k_3}e_{j_2k_3}\right)\\
	&=&O\left(\frac{r}{nT}\right)
	+\frac{\sqrt{C_1}}{n^2 T^4}\sum_{s_1\neq t}^{T}\sum_{s_2\neq s_1,t}^{T}\sum_{t=1}^{T}\sum_{j_1\neq i_1}^{n}\sum_{i_1=1}^{n}\sum_{j_2\neq i_1}^{n}\left(\sum_{k_2=1}^{r}h_{k_2}e_{i_1k_2}e_{j_1k_2}\right)\left(\sum_{k_3=1}^{r}h_{k_3}e_{i_1k_3}e_{j_2k_3}\right)\\
	&&+\frac{\sqrt{C_1}}{n^2 T^4}\sum_{s_1\neq t}^{T}\sum_{s_2\neq s_1,t}^{T}\sum_{t=1}^{T}\sum_{j_1\neq i_1}^{n}\sum_{i_1=1}^{n}\sum_{j_2\neq i_2}^{n}\sum_{i_2=1}^{n}\left(\sum_{k_1=1}^{r}h_{k_1}e_{i_1k_1}e_{i_2k_1}\right)\left(\sum_{k_2=1}^{r}h_{k_2}e_{i_1k_2}e_{j_1k_2}\right)\left(\sum_{k_3=1}^{r}h_{k_3}e_{i_2k_3}e_{j_2k_3}\right)\\
	&=&O\left(\frac{r}{nT}\right)+O\left(\frac{r}{nT}\right)+\frac{\sqrt{C_1}}{n^2 T^4}\sum_{s_1\neq t}^{T}\sum_{s_2\neq s_1,t}^{T}\sum_{t=1}^{T}\sum_{j_1\neq i_1}^{n}\sum_{i_1=1}^{n}\sum_{j_2\neq i_2}^{n}\sum_{i_2=1}^{n}\left(\sum_{k=1}^{r}h_k^3 e_{i_1k}^2 e_{i_2k}^2 e_{j_1k}e_{j_2k}\right)\\
	&=&O\left(\frac{r}{nT}\right)+O\left(\frac{r}{nT}\right)+O\left(\frac{r}{nT}\right)=O\left(\frac{r}{nT}\right).
\end{eqnarray*}
Therefore, the first term of (\ref{f.1}) is
\begin{eqnarray*}
	\frac{1}{n T^2}\sum_{s=1}^{T}\sum_{t=1}^{T}v_t^{\prime}v_s v_s^{\prime}\tilde{x}_t
	&=&\frac{1}{n T^2}\sum_{t=1}^{T}\sum_{i=1}^{n}v_{it}^3\tilde{x}_{it}^{\prime}+\frac{1}{n T^2}\sum_{t=1}^{T}\sum_{j\neq i}^{n}\sum_{i=1}^{n}v_{it}^2 v_{jt}\tilde{x}_{jt}^{\prime}\\
	&& +\frac{1}{n T^2}\sum_{s\neq t}^{T}\sum_{t=1}^{T}\sum_{i=1}^{n}v_{it}v_{is}^2\tilde{x}_{it}^{\prime}+\frac{1}{n T^2}\sum_{s\neq t}^{T}\sum_{t=1}^{T}\sum_{j\neq i}^{n}\sum_{i=1}^{n}v_{it}v_{is}v_{js}\tilde{x}_{jt}^{\prime}\\
	&=&O_p\left(\frac{1}{T}\right)+\left(O_p\left(\frac{r}{T\sqrt{n}}\right)+O_p\left(\frac{n}{T\sqrt{T}}\right)\right)+O_p\left(\frac{1}{\sqrt{T}}\right)+O_p\left(\sqrt{\frac{r}{nT}}\right)\\
	&=&O_p\left(\frac{r}{T\sqrt{n}}\right)+O_p\left(\frac{n}{T\sqrt{T}}\right)+O_p\left(\frac{1}{\sqrt{T}}\right).
\end{eqnarray*}
Consider the second term of (\ref{f.1}),
\begin{eqnarray*}
	\frac{1}{n T^2}\sum_{s=1}^{T}\sum_{t=1}^{T}\bar{v}^{\prime}v_s v_s^{\prime}\tilde{x}_t=O_p\left(\frac{r}{T\sqrt{n}}\right)+O_p\left(\frac{n}{T\sqrt{T}}\right)+O_p\left(\frac{1}{\sqrt{T}}\right).
\end{eqnarray*}
The calculation of the second term is similar to that of the first term. Thus,
\begin{eqnarray*}
	\frac{1}{n}tr(SA_1)
	&=&\left(O_p\left(\frac{r}{T\sqrt{n}}\right)+O_p\left(\frac{n}{T\sqrt{T}}\right)+O_p\left(\frac{1}{\sqrt{T}}\right)\right)O_p\left(\frac{1}{\sqrt{nT}}\right)\\
	&=&O_p\left(\frac{r}{nT\sqrt{T}}\right)+O_p\left(\frac{\sqrt{n}}{T^2}\right)+O_p\left(\frac{1}{T\sqrt{n}}\right).
\end{eqnarray*}

\subsubsection*{Proof of part (b)}
\begin{eqnarray*}
	\frac{1}{n}tr(SA_3)=\frac{1}{nT}\left(\tilde{\beta}-\beta\right)^{\prime}\left(\sum_{t=1}^{T}\tilde{x}_t^{\prime}S\tilde{x}_t\right)\left(\tilde{\beta}-\beta\right),
\end{eqnarray*}
where
\begin{eqnarray*}
	\frac{1}{nT}\sum_{t=1}^{T}\tilde{x}_t^{\prime}S\tilde{x}_t&=&\frac{1}{nT}\sum_{t=1}^{T}\tilde{x}_t^{\prime}\left(\frac{1}{T}\sum_{s=1}^{T}v_s v_s^{\prime}\right)\tilde{x}_t\\
	&=&\frac{1}{nT^2}\sum_{t=1}^{T}\tilde{x}_t^{\prime}v_t v_t^{\prime}\tilde{x}_t+\frac{1}{nT^2}\sum_{s\neq t}^{T}\sum_{t=1}^{T}\tilde{x}_t^{\prime}v_s v_s^{\prime}\tilde{x}_t\\
	&=&\frac{1}{nT^2}\sum_{t=1}^{T}\sum_{i=1}^{n}\tilde{x}_{it}^2v_{it}^2+\frac{1}{nT^2}\sum_{t=1}^{T}\sum_{j\neq i}^{n}\sum_{i=1}^{n}\tilde{x}_{it}v_{it}v_{jt}\tilde{x}_{jt}\\
	&&+\frac{1}{nT^2}\sum_{s\neq t}^{T}\sum_{t=1}^{T}\sum_{i=1}^{n}\tilde{x}_{it}^2v_{is}^2+\frac{1}{nT^2}\sum_{s\neq t}^{T}\sum_{t=1}^{T}\sum_{j\neq i}^{n}\sum_{i=1}^{n}\tilde{x}_{it}v_{is}v_{js}\tilde{x}_{jt}\\
	&=&O_p\left(\frac{1}{T}\right)+O_p\left(\frac{n}{T}\right)+O_p\left(1\right)+O_p\left(r\right).
\end{eqnarray*}

$\frac{1}{nT^2}\sum_{s\neq t}^{T}\sum_{t=1}^{T}\sum_{j\neq i}^{n}\sum_{i=1}^{n}\tilde{x}_{it}v_{is}v_{js}\tilde{x}_{jt}=O_p\left(r\right)$, because
\begin{eqnarray*}
	&&E\left(\frac{1}{nT^2}\sum_{s\neq t}^{T}\sum_{t=1}^{T}\sum_{j\neq i}^{n}\sum_{i=1}^{n}\tilde{x}_{it}v_{is}v_{js}\tilde{x}_{jt}\right)^2\\
	&=&\frac{1}{n^2T^4}\sum_{s_1\neq t_1}^{T}\sum_{t_1=1}^{T}\sum_{j_1\neq i_1}^{n}\sum_{i_1=1}^{n}\sum_{s_2\neq t_2}^{T}\sum_{t_2=1}^{T}\sum_{j_2\neq i_2}^{n}\sum_{i_2=1}^{n}E\tilde{x}_{i_1t_1}v_{i_1s_1}v_{j_1s_1}\tilde{x}_{j_1t_1}\tilde{x}_{i_2t_2}v_{i_2s_2}v_{j_2s_2}\tilde{x}_{j_2t_2}\\
	&=&\frac{1}{n^2T^4}\sum_{s_1\neq t_1}^{T}\sum_{t_1=1}^{T}\sum_{j_1\neq i_1}^{n}\sum_{i_1=1}^{n}\sum_{s_2\neq t_2}^{T}\sum_{t_2=1}^{T}\sum_{j_2\neq i_2}^{n}\sum_{i_2=1}^{n}Ev_{i_1s_1}v_{j_1s_1}v_{i_2s_2}v_{j_2s_2}E\tilde{x}_{i_1t_1}\tilde{x}_{j_1t_1}\tilde{x}_{i_2t_2}\tilde{x}_{j_2t_2}\\
	&\leq&\frac{C_1}{n^2T^4}\sum_{s_1\neq t_1}^{T}\sum_{t_1=1}^{T}\sum_{j_1\neq i_1}^{n}\sum_{i_1=1}^{n}\sum_{s_2\neq t_2}^{T}\sum_{t_2=1}^{T}\sum_{j_2\neq i_2}^{n}\sum_{i_2=1}^{n}Ev_{i_1s_1}v_{j_1s_1}v_{i_2s_2}v_{j_2s_2}\\
	&=&\frac{C_1(T-1)^2}{n^2T^4}\sum_{s_1=1}^{T}\sum_{j_1\neq i_1}^{n}\sum_{i_1=1}^{n}\sum_{s_2=1}^{T}\sum_{j_2\neq i_2}^{n}\sum_{i_2=1}^{n}Ev_{i_1s_1}v_{j_1s_1}v_{i_2s_2}v_{j_2s_2}\\
	&=&\frac{C_1(T-1)^2}{n^2T^4}\sum_{s_1=1}^{T}\sum_{j_1\neq i_1}^{n}\sum_{i_1=1}^{n}\sum_{j_2\neq i_2}^{n}\sum_{i_2=1}^{n}Ev_{i_1s_1}v_{j_1s_1}v_{i_2s_1}v_{j_2s_1}\\
	&&+\frac{C_1(T-1)^2}{n^2T^4}\sum_{s_1=1}^{T}\sum_{j_1\neq i_1}^{n}\sum_{i_1=1}^{n}\sum_{s_2\neq s_1}^{T}\sum_{j_2\neq i_2}^{n}\sum_{i_2=1}^{n}Ev_{i_1s_1}v_{j_1s_1}Ev_{i_2s_2}v_{j_2s_2}\\
	&=&O\left(\frac{r^2}{T}\right)+\frac{C_1T(T-1)^3}{n^2T^4}\sum_{j_1\neq i_1}^{n}\sum_{i_1=1}^{n}\sum_{j_2\neq i_2}^{n}\sum_{i_2=1}^{n}\Sigma_{i_1j_1}\Sigma_{i_2j_2}\\
	&=&O\left(\frac{r^2}{T}\right)+\frac{C_1T(T-1)^3}{n^2T^4}\left(\sum_{j\neq i}^{n}\sum_{i=1}^{n}\Sigma_{ij}\right)^2\\
	&=&O\left(\frac{r^2}{T}\right)+\frac{C_1T(T-1)^3}{n^2T^4}\left(\sum_{j\neq i}^{n}\sum_{i=1}^{n}\left(\sum_{k=1}^{r}h_ke_{ik}e_{jk}\right)\right)^2\\
	&=&O\left(\frac{r^2}{T}\right)+\frac{C_1T(T-1)^3}{n^2T^4}\left(O\left(nr\right)\right)^2\\
	&=&O(r^2).
\end{eqnarray*}

Here we apply Cauchy-Schwarz Inequality, $(\sum_{i=1}^{n}e_{ik})^2\leq n\sum_{i=1}^{n}e_{ik}^2=n$.

\subsubsection*{Proof of part (g)}
\begin{eqnarray*}
	\frac{1}{n}tr(SA_0)&=&\frac{1}{n}tr(S\bar{v}\bar{v}^{\prime})\\
	&=&\frac{1}{n}tr\left((\frac{1}{T}\sum_{t=1}^{T}v_tv_t^{\prime})(\frac{1}{T}\sum_{s=1}^{T}v_s)(\frac{1}{T}\sum_{\tau=1}^{T}v_\tau)^{\prime}\right)\\
	&=&\frac{1}{nT^3}v_t^{\prime}v_sv_{\tau}^{\prime}v_t\\
	&=&\frac{1}{nT^3}\sum_{t=1}^{T}\sum_{s=1}^{T}\sum_{\tau=1}^{T}\sum_{i=1}^{n}\sum_{j=1}^{n}v_{it}v_{is}v_{j\tau}v_{jt}\\
	&=&\frac{1}{nT^3}\sum_{t=1}^{T}\sum_{s=1}^{T}\sum_{i=1}^{n}\sum_{j=1}^{n}v_{it}v_{is}v_{js}v_{jt}+\frac{1}{nT^3}\sum_{t=1}^{T}\sum_{s=1}^{T}\sum_{\tau\neq s}^{T}\sum_{i=1}^{n}\sum_{j=1}^{n}v_{it}v_{is}v_{j\tau}v_{jt}\\
	&=&\frac{1}{nT}trS^2\\&&+\frac{1}{nT^3}\left(\sum_{t=1}^{T}\sum_{\tau\neq t}^{T}\sum_{i=1}^{n}\sum_{j=1}^{n}v_{it}^2v_{j\tau}v_{jt}+\sum_{t=1}^{T}\sum_{s\neq t}^{T}\sum_{i=1}^{n}\sum_{j=1}^{n}v_{jt}^2v_{it}v_{is}+\sum_{t=1}^{T}\sum_{s\neq t}^{T}\sum_{\tau\neq s,t}^{T}\sum_{i=1}^{n}\sum_{j=1}^{n}v_{it}v_{is}v_{j\tau}v_{jt}\right)\\
	&=&\frac{1}{nT}trS^2+R_1.
\end{eqnarray*}

\subsubsection*{Proof of part (h)}
\begin{eqnarray*}
	\frac{1}{n}tr(A_0^2)&=&\frac{1}{n}\bar{v}^{\prime}\bar{v}\bar{v}^{\prime}\bar{v}=n\left(\frac{1}{n}\bar{v}^{\prime}\bar{v}\right)^2\\
	&=&n\left(\frac{1}{nT^2}\sum_{i=1}^{n}\sum_{t=1}^{T}v_{it}^2+\frac{1}{nT^2}\sum_{i=1}^{n}\sum_{s\neq t}^{T}\sum_{t=1}^{T}v_{is}v_{it}\right)^2\\
	&=&n\left(\frac{1}{nT}trS+O_p\left(\frac{1}{T\sqrt{n}}\right)\right)^2\\
	&=&\frac{n}{T^2}\left(\frac{1}{n}trS\right)^2+O_p\left(\frac{\sqrt{n}}{T^2}\right).
\end{eqnarray*}
$\frac{1}{nT^2}\sum_{i=1}^{n}\sum_{s\neq t}^{T}\sum_{t=1}^{T}v_{is}v_{it}=O_p\left(\frac{1}{T\sqrt{n}}\right)$, because
\begin{eqnarray*}
	E\left(\frac{1}{nT^2}\sum_{i=1}^{n}\sum_{s\neq t}^{T}\sum_{t=1}^{T}v_{is}v_{it}\right)^2
	&=&\frac{1}{n^2T^4}\sum_{i=1}^{n}\sum_{s\neq t}^{T}\sum_{t=1}^{T}\sum_{j=1}^{n}\sum_{s_2\neq t_2}^{T}\sum_{t_2=1}^{T}Ev_{is}v_{it}v_{js_2}v_{jt_2}\\
	&=&\frac{1}{n^2T^4}\sum_{i=1}^{n}\sum_{j=1}^{n}\sum_{s\neq t}^{T}\sum_{t=1}^{T}Ev_{is}v_{it}v_{js}v_{jt}\\
	&=&\frac{1}{n^2T^4}\sum_{i=1}^{n}\sum_{j=1}^{n}\sum_{s\neq t}^{T}\sum_{t=1}^{T}\Sigma_{ij}^2\\
	&=&\frac{T-1}{n^2T^3}\sum_{i=1}^{n}\Sigma_{ii}^2+\frac{T-1}{n^2T^3}\sum_{i=1}^{n}\sum_{j\neq i}^{n}\Sigma_{ij}^2\\
	&=&O_p\left(\frac{1}{nT^2}\right)+\frac{T-1}{n^2T^3}\sum_{i=1}^{n}\sum_{j\neq i}^{n}\left(\sum_{k=1}^{r}h_ke_{ik}e_{jk}\right)^2\\
	&\leq&O_p\left(\frac{1}{nT^2}\right)+\frac{T-1}{n^2T^3}\sum_{i=1}^{n}\sum_{j=1}^{n}\left(\sum_{k=1}^{r}h_ke_{ik}e_{jk}\right)^2\\
	&=&O_p\left(\frac{1}{nT^2}\right)+\frac{T-1}{n^2T^3}\sum_{i=1}^{n}\sum_{j=1}^{n}\sum_{k=1}^{r}h_k^2 e_{ik}^2 e_{jk}^2\\
	&=&O_p(\frac{1}{nT^2})+\frac{T-1}{n^2T^3}\sum_{k=1}^{r}h_k^2\\
	&=&O_p\left(\frac{1}{nT^2}\right)+O_p\left(\frac{r}{n^2T^2}\right).
\end{eqnarray*}
The result comes from the orthogonality of $\textbf{e}_j,j=1,\ldots,r$ and $\Vert\textbf{e}_j\Vert=1$.

\begin{lemma}\label{S_diverge_factor}
	Under Assumptions 1, 2 and the divergent number of weak factors alternative,
	
	(a)$\frac{1}{n}tr(S)=\sigma^2(1+\frac{1}{n}\sum_{j=1}^{r}h_j)+O_p(\frac{1}{\sqrt{T}})=O_p(1)$;
	
	(b)$\frac{1}{n}tr(S^2)=O_p\left(\frac{n}{T}\right)+O_p\left(1\right)$;
	
	(c) $W_1=\frac{1}{n}tr(\hat{S})-\frac{1}{n}tr(S)=-\frac{1}{nT}trS+R_2+O_p(\frac{1}{nT})$;
	
	(d) $W_2=\frac{1}{n}tr(\hat{S}^2)-\frac{1}{n}tr(S^2)=-2\frac{1}{n}tr(SA_0)+\frac{1}{n}tr(A_0^2)+O_p(1)$.
\end{lemma}

\begin{lemma}\label{Remainder}
	Under Assumptions 1,2 and the divergent number of weak factors alternative,
	\begin{eqnarray*}
		R_1&=&O_p\left(\frac{n}{T^2}\right),\\
		R_2&=&O_p\left(\frac{1}{T\sqrt{n}}\right),
	\end{eqnarray*}
\end{lemma}

\subsubsection*{Proof of Lemma \ref{Remainder}}
\begin{eqnarray*}
	&&E\left(\frac{1}{nT^3}\sum_{t=1}^{T}\sum_{s\neq t}^{T}\sum_{i=1}^{n}\sum_{j=1}^{n}v_{jt}^2v_{it}v_{is}\right)^2\\
	&=&\frac{1}{n^2T^6}\sum_{t_1=1}^{T}\sum_{t_2=1}^{T}\sum_{s\neq t_1,t_2}^{T}\sum_{i_1=1}^{n}\sum_{j_1=1}^{n}\sum_{i_2=1}^{n}\sum_{j_2=1}^{n}Ev_{j_1t_1}^2v_{i_1t_1}v_{j_2t_2}^2v_{i_2t_2}Ev_{i_1s}v_{i_2s}\\
	&=&\frac{1}{n^2T^6}\sum_{t_1=1}^{T}\sum_{s\neq t_1}^{T}\sum_{i_1=1}^{n}\sum_{j_1=1}^{n}\sum_{i_2=1}^{n}\sum_{j_2=1}^{n}Ev_{j_1t_1}^2v_{i_1t_1}v_{j_2t_1}^2v_{i_2t_1}Ev_{i_1s}v_{i_2s}\\
	&&+\frac{1}{n^2T^6}\sum_{t_1=1}^{T}\sum_{t_2\neq t_1}^{T}\sum_{s\neq t_1,t_2}^{T}\sum_{i_1=1}^{n}\sum_{j_1=1}^{n}\sum_{i_2=1}^{n}\sum_{j_2=1}^{n}Ev_{j_1t_1}^2v_{i_1t_1}Ev_{j_2t_2}^2v_{i_2t_2}Ev_{i_1s}v_{i_2s}\\
	&=&O\left(\frac{n^2}{T^4}\right)+\frac{1}{n^2T^6}\sum_{t_1=1}^{T}\sum_{t_2\neq t_1}^{T}\sum_{s\neq t_1,t_2}^{T}\sum_{i_1=1}^{n}\sum_{j_1=1}^{n}\sum_{i_2=1}^{n}\sum_{j_2=1}^{n}Ev_{j_1t_1}^2v_{i_1t_1}Ev_{j_2t_2}^2v_{i_2t_2}Ev_{i_1s}v_{i_2s}.
\end{eqnarray*}

Plug in $Ev_{j_1t_1}^2v_{i_1t_1}=\sum_{k=1}^{r}\gamma_{j_1k}^2\gamma_{i_1k}Ef_{t_1k}^3$ and $Ev_{i_1s}v_{i_2s}=\sum_{k=1}^{r}\gamma_{i_1k}\gamma_{i_2k}\sigma_k^2$,

\begin{eqnarray*}
	&&E\left(\frac{1}{nT^3}\sum_{t=1}^{T}\sum_{s\neq t}^{T}\sum_{i=1}^{n}\sum_{j=1}^{n}v_{jt}^2v_{it}v_{is}\right)^2\\
	&=&O\left(\frac{n^2}{T^4}\right)\\&&+O\left(\frac{1}{n^2T^3}\right)\sum_{i_1=1}^{n}\sum_{j_1=1}^{n}\sum_{i_2=1}^{n}\sum_{j_2=1}^{n}\left(\sum_{k_1=1}^{r}\gamma_{j_1k_1}^2\gamma_{i_1k_1}Ef_{t_1k_1}^3\right)\left(\sum_{k_2=1}^{r}\gamma_{j_2k_2}^2\gamma_{i_2k_2}Ef_{t_1k_2}^3\right)\left(\sum_{k_3=1}^{r}\gamma_{i_1k_3}\gamma_{i_2k_3}\sigma_{k_3}^2\right)\\
	&=&O\left(\frac{n^2}{T^4}\right)+O\left(\frac{1}{n^2T^3}\right)\sum_{i_1=1}^{n}\sum_{j_1=1}^{n}\sum_{i_2=1}^{n}\sum_{j_2=1}^{n}\sum_{k=1}^{r}\gamma_{j_1k}^2\gamma_{i_1k}^2\gamma_{j_2k}^2\gamma_{i_2k}^2(Ef_{t_1k}^3)^2 \sigma_{k}^2\\
	&=&O\left(\frac{n^2}{T^4}\right)+O\left(\frac{r}{n^2T^3}\right)=O\left(\frac{n^2}{T^4}\right).
\end{eqnarray*}

\begin{eqnarray*}
	&&E\left(\frac{1}{nT^3}\sum_{t=1}^{T}\sum_{s\neq t}^{T}\sum_{\tau\neq s,t}^{T}\sum_{i=1}^{n}\sum_{j=1}^{n}v_{it}v_{is}v_{j\tau}v_{jt}\right)^2\\
	&=&\frac{1}{n^2T^6}\sum_{t_1=1}^{T}\sum_{t_2=1}^{T}\sum_{s\neq t_1,t_2}^{T}\sum_{\tau\neq s,t_1,t_2}^{T}\sum_{i_1=1}^{n}\sum_{j_1=1}^{n}\sum_{i_2=1}^{n}\sum_{j_2=1}^{n}Ev_{i_1t_1}v_{j_1t_1}v_{i_2t_2}v_{j_2t_2}Ev_{i_1s}v_{i_2s}Ev_{j_1\tau}v_{j_2\tau}\\
	&=&\frac{1}{n^2T^6}\sum_{t_1=1}^{T}\sum_{s\neq t_1}^{T}\sum_{\tau\neq s,t_1}^{T}\sum_{i_1=1}^{n}\sum_{j_1=1}^{n}\sum_{i_2=1}^{n}\sum_{j_2=1}^{n}Ev_{i_1t_1}v_{j_1t_1}v_{i_2t_1}v_{j_2t_1}Ev_{i_1s}v_{i_2s}Ev_{j_1\tau}v_{j_2\tau}\\
	&&+\frac{1}{n^2T^6}\sum_{t_1=1}^{T}\sum_{t_2\neq t_1}^{T}\sum_{s\neq t_1,t_2}^{T}\sum_{\tau\neq s,t_1,t_2}^{T}\sum_{i_1=1}^{n}\sum_{j_1=1}^{n}\sum_{i_2=1}^{n}\sum_{j_2=1}^{n}Ev_{i_1t_1}v_{j_1t_1}Ev_{i_2t_2}v_{j_2t_2}Ev_{i_1s}v_{i_2s}Ev_{j_1\tau}v_{j_2\tau}.
\end{eqnarray*}

Plug in $Ev_{it}v_{jt}=\sum_{k=1}^{r}\gamma_{ik}\gamma_{jk}\sigma_k^2$ and $Ev_{i_1t_1}v_{j_1t_1}v_{i_2t_1}v_{j_2t_1}=\sum_{k=1}^{r}\gamma_{i_1k}\gamma_{j_1k}\gamma_{i_2k}\gamma_{j_2k}Ef_{t_1k}^4+3\sum_{k_2\neq k_1}^{r}\sum_{k_1=1}^{r}\gamma_{i_1k_1}\gamma_{i_2k_1}\gamma_{j_1k_2}\gamma_{j_2k_2}\sigma_{k_1}^2\sigma_{k_2}^2$,
\begin{eqnarray*}
	&&E\left(\frac{1}{nT^3}\sum_{t=1}^{T}\sum_{s\neq t}^{T}\sum_{\tau\neq s,t}^{T}\sum_{i=1}^{n}\sum_{j=1}^{n}v_{it}v_{is}v_{j\tau}v_{jt}\right)^2\\
	&=&O\left(\frac{1}{n^2T^3}\right)\sum_{i_1=1}^{n}\sum_{j_1=1}^{n}\sum_{i_2=1}^{n}\sum_{j_2=1}^{n}\left(\sum_{k=1}^{r}\gamma_{i_1k}\gamma_{j_1k}\gamma_{i_2k}\gamma_{j_2k}Ef_{t_1k}^4+3\sum_{k_2\neq k_1}^{r}\sum_{k_1=1}^{r}\gamma_{i_1k_1}\gamma_{i_2k_1}\gamma_{j_1k_2}\gamma_{j_2k_2}\sigma_{k_1}^2\sigma_{k_2}^2\right)\\
	&&\left(\sum_{k=1}^{r}\gamma_{i_1k}\gamma_{i_2k}\sigma_k^2\right)\left(\sum_{k=1}^{r}\gamma_{j_1k}\gamma_{j_2k}\sigma_k^2\right)\\
	&&+O\left(\frac{1}{n^2T^2}\right)\sum_{i_1=1}^{n}\sum_{j_1=1}^{n}\sum_{i_2=1}^{n}\sum_{j_2=1}^{n}\sum_{k=1}^{r}\gamma_{i_1k}^2\gamma_{i_2k}^2\gamma_{j_1k}^2\gamma_{j_2k}^2\\
	&=&O\left(\frac{1}{n^2T^3}\right)\sum_{i_1=1}^{n}\sum_{j_1=1}^{n}\sum_{i_2=1}^{n}\sum_{j_2=1}^{n}\left(\sum_{k=1}^{r}\gamma_{i_1k}^2\gamma_{j_1k}^2\gamma_{i_2k}^2\gamma_{j_2k}^2+\sum_{k_2\neq k_1}^{r}\sum_{k_1=1}^{r}\gamma_{i_1k_1}^2\gamma_{i_2k_1}^2\gamma_{j_1k_2}^2\gamma_{j_2k_2}^2\right)+O\left(\frac{r}{n^2T^2}\right)\\
	&=&O\left(\frac{r^2}{n^2T^3}\right)+O\left(\frac{r}{n^2T^2}\right).
\end{eqnarray*}

\begin{eqnarray*}
	ER_2^2&=&E\left(-\frac{1}{nT^2}\sum_{t=1}^{T}\sum_{s\neq t}^{T}\sum_{i=1}^{n}v_{is}v_{it}\right)^2\\
	&=&\frac{1}{n^2T^4}\sum_{t=1}^{T}\sum_{s\neq t}^{T}\sum_{i_1=1}^{n}Ev_{i_1t}^2Ev_{i_1s}^2+\frac{1}{n^2T^4}\sum_{t=1}^{T}\sum_{s\neq t}^{T}\sum_{i_1=1}^{n}\sum_{i_2\neq i_1}^{n}Ev_{i_1t}v_{i_2t}Ev_{i_1s}v_{i_2s}\\
	&=&O\left(\frac{1}{nT^2}\right)+\frac{T(T-1)}{n^2T^4}\sum_{i_1=1}^{n}\sum_{i_2\neq i_1}^{n}(\Sigma_{i_1,i_2})^2\\
	&\leq&O\left(\frac{1}{nT^2}\right)+\frac{T(T-1)}{n^2T^4}\sum_{i_1=1}^{n}\sum_{i_2=1}^{n}(\Sigma_{i_1,i_2})^2\\
	&=&O\left(\frac{1}{nT^2}\right)+\frac{T(T-1)}{n^2T^4}\sum_{i_1=1}^{n}\sum_{i_2=1}^{n}\left(\sum_{k=1}^{r}h_ke_{i_1k}e_{i_2k}\right)^2\\
	&=&O\left(\frac{1}{nT^2}\right)+\frac{T(T-1)}{n^2T^4}\sum_{k=1}^{r}h_k^2\\
	&=&O\left(\frac{1}{nT^2}\right)+O\left(\frac{r}{n^2T^2}\right).
\end{eqnarray*}

\begin{eqnarray*}
	&&E\left(\frac{1}{nT^3}\sum_{t=1}^{T}\sum_{s\neq t}^{T}\sum_{i=1}^{n}\sum_{j=1}^{n}v_{jt}^2v_{it}v_{is}\right)^2\\
	&=&\frac{1}{n^2T^6}\sum_{t_1=1}^{T}\sum_{t_2=1}^{T}\sum_{s\neq t_1,t_2}^{T}\sum_{i_1=1}^{n}\sum_{j_1=1}^{n}\sum_{i_2=1}^{n}\sum_{j_2=1}^{n}Ev_{j_1t_1}^2v_{i_1t_1}v_{i_1s}v_{j_2t_2}^2v_{i_2t_2}v_{i_2s}\\
	&=&\frac{1}{n^2T^6}\sum_{t_1=1}^{T}\sum_{t_2=1}^{T}\sum_{s\neq t_1,t_2}^{T}\sum_{i_1=1}^{n}\sum_{j_1=1}^{n}\sum_{i_2=1}^{n}\sum_{j_2=1}^{n}Ev_{j_1t_1}^2v_{i_1t_1}v_{j_2t_2}^2v_{i_2t_2}Ev_{i_1s}v_{i_2s}\\
	&=&\frac{1}{n^2T^6}\sum_{t_1=1}^{T}\sum_{s\neq t_1,t_2}^{T}\sum_{i_1=1}^{n}\sum_{j_1=1}^{n}\sum_{i_2=1}^{n}\sum_{j_2=1}^{n}Ev_{j_1t_1}^2v_{i_1t_1}v_{j_2t_1}^2v_{i_2t_1}Ev_{i_1s}v_{i_2s}\\&&+\frac{1}{n^2T^6}\sum_{t_1=1}^{T}\sum_{t_2\neq t_1}^{T}\sum_{s\neq t_1,t_2}^{T}\sum_{i_1=1}^{n}\sum_{j_1=1}^{n}\sum_{i_2=1}^{n}\sum_{j_2=1}^{n}Ev_{j_1t_1}^2v_{i_1t_1}Ev_{j_2t_2}^2v_{i_2t_2}Ev_{i_1s}v_{i_2s}.
\end{eqnarray*}
Finally, we have
\begin{eqnarray*}
	&&T(\hat{U}-U)\\
	&=&\frac{TW_2(\frac{1}{n}trS)^2-2TW_1\frac{1}{n}trS\frac{1}{n}tr S^2-T W_1^2\frac{1}{n}tr S^2}{\left(\frac{1}{n}trS+W_1\right)^2\left(\frac{1}{n}trS\right)^2}\\
	&=&\frac{T\frac{1}{n}trS\left(W_2\frac{1}{n}trS-2W_1\frac{1}{n}trS^2\right)-TW_1^2\frac{1}{n}trS^2}{\left(\frac{1}{n}trS+O_p\left(\frac{1}{T}\right)\right)^2\left(\frac{1}{n}trS\right)^2}\\
	&=&\frac{T\frac{1}{n}trS\left[\left(-2\frac{1}{n}tr(SA_0)+\frac{1}{n}tr(A_0^2)+O_p(\frac{r}{nT})\right)\frac{1}{n}trS-2\left(-\frac{1}{nT}trS+R_2+O_p(\frac{1}{nT})\right)\frac{1}{n}trS^2\right]+O_p\left(\frac{1}{T}\right)}{\left(\frac{1}{n}trS+O_p\left(\frac{1}{T}\right)\right)^2\left(\frac{1}{n}trS\right)^2}\\
	&=&\frac{T\frac{1}{n}trS\left[\left(-\frac{2}{nT}tr(S^2)+R_1\right)\frac{1}{n}trS-2\left(-\frac{1}{nT}trS+R_2\right)\frac{1}{n}trS^2\right]+O_p\left(\frac{1}{T}\right)+O_p\left(\frac{r}{n}\right)}{\left(\frac{1}{n}trS+O_p\left(\frac{1}{T}\right)\right)^2\left(\frac{1}{n}trS\right)^2}+\frac{T\frac{1}{n}tr(A_0^2)}{\left(\frac{1}{n}trS+O_p\left(\frac{1}{T}\right)\right)^2}\\
	&=&\frac{T\frac{1}{n}trS\left[R_1\frac{1}{n}trS-2R_2\frac{1}{n}trS^2\right]+O_p\left(\frac{r}{n}\right)+O_p\left(\frac{1}{T}\right)}{\left(\frac{1}{n}trS+O_p\left(\frac{1}{T}\right)\right)^2\left(\frac{1}{n}trS\right)^2}+\frac{T\frac{1}{n}tr(A_0^2)}{\left(\frac{1}{n}trS+O_p\left(\frac{1}{T}\right)\right)^2}\\
	&=&\frac{T\left[R_1\frac{1}{n}trS-2R_2\frac{1}{n}trS^2\right]}{\left(\frac{1}{n}trS+O_p\left(\frac{1}{T}\right)\right)^2\left(\frac{1}{n}trS\right)}+\frac{T\frac{1}{n}tr(A_0^2)}{\left(\frac{1}{n}trS+O_p\left(\frac{1}{T}\right)\right)^2}+O_p\left(\frac{r}{n}\right)+O_p\left(\frac{1}{T}\right)\\
	&=&\frac{T\left[R_1\frac{1}{n}trS-2R_2\frac{1}{n}trS^2\right]}{\left(\frac{1}{n}trS+O_p\left(\frac{1}{T}\right)\right)^2\left(\frac{1}{n}trS\right)}+\frac{T\left(\frac{n}{T^2}\left(\frac{1}{n}trS\right)^2+O_p\left(\frac{\sqrt{n}}{T^2}\right)\right)}{\left(\frac{1}{n}trS+O_p\left(\frac{1}{T}\right)\right)^2}+O_p\left(\frac{r}{n}\right)+O_p\left(\frac{1}{T}\right)\\
	&=&\frac{T\left[R_1\frac{1}{n}trS-2R_2\frac{1}{n}trS^2\right]}{\left(\frac{1}{n}trS+O_p\left(\frac{1}{T}\right)\right)^2\left(\frac{1}{n}trS\right)}+c_T+O_p\left(\frac{\sqrt{n}}{T}\right)+O_p\left(\frac{r}{n}\right)+O_p\left(\frac{1}{T}\right)\\
	&=&O_p\left(\frac{n}{T}\right)+O_p\left(\frac{1}{\sqrt{n}}\right)+c_T+O_p\left(\frac{r}{n}\right)+O_p\left(\frac{\sqrt{n}}{T}\right).
\end{eqnarray*}
 The proof of Proposition \ref{diff_between_U_diverge_factor} is complete.

\section{Proof of Proposition \ref{diff_between_U_not_weak}}

\begin{proposition}\label{S}
	Suppose Assumptions 1-4 hold. Under the LPA scheme, $c_T=n/T\rightarrow c$. Under a neither strong nor weak factor model alternative, that is $h_j/n^{\alpha}\rightarrow d_j,0<\alpha<1,0<d_j<\infty$, we have
	\begin{eqnarray*}
		&& \frac{1}{n}trS=\sigma^2+O_p\left(\frac{1}{n^{1-\alpha}}\right)+O_p\left(\frac{1}{\sqrt{T}}\right),\\
		&& \frac{1}{n}trS^2=\left(\frac{n}{T}+1\right)\sigma^4+O_p\left(n^{2\alpha-1}\right)+O_p\left(\frac{1}{\sqrt{T}}\right),\\
		&& \frac{1}{n}tr\hat{S}-\frac{1}{n}trS=-\frac{\sigma^2}{T}+O_p\left(\frac{n^{\alpha-1}}{T}\right),\\
		&& \frac{1}{n}tr\hat{S}^2-\frac{1}{n}trS^2=-\frac{2}{T}\sigma^4-\frac{n}{T^2}\sigma^4+O_p\left(\frac{n^{2\alpha-1}}{T}\right)+O_p\left(\frac{1}{T\sqrt{T}}\right).
	\end{eqnarray*}
\end{proposition}
Then the proof is complete.

\end{document}


\newcommand{\red}{\color{DarkRed}}
\newcommand\cov{\mathop{\text{cov}}}
\newcommand\var{\mathop{\text{var}}}
\newcommand\diag{\mathop{\text{diag}}}
\newcommand\cvd{\xrightarrow{~d~}{}}
\renewcommand{\(}{\left(}
\renewcommand{\)}{\right)}

\begin{center}
	{\bf \large On John's test for sphericity in large panel data models\\ (Supplementary Material)}
	\\[4mm]
	\title
	\author{Zhaoyuan LI\\[1mm] The Chinese University of Hong Kong, Shenzhen}\\[2mm] July 3, 2022
\end{center}

\bigskip

\section{Asymptotic power under general form of covariance}

Based on the power analysis of the GRJ test under the weak factor alternative in Section 4.1, i.e., the disturbances follow a factor model with fixed number of weak factors, we continue to explore the power behavior of the GRJ test under general form of covariance with bouned spectral norm $\Sigma_n\neq \sigma_{\nu}^2 \mathbf{I}_n$ which covers the weak factor alternative as a special case. We first develop the limiting distribution of $U$ for general covariance matrix and investigate the asymptotic power of John's test for raw data $\{\nu_{it} \}$. We will write $a_n\sim b_n$ for positive numbers $a_n$ and $b_n$ if $a_n/b_n\to 1$ as $n\to \infty$. 

\begin{lemma}\label{high_alter}
	Suppose the moments of $H_n$, the ESD of covariance matrix ${\Sigma}_n$, converge to the corresponding moments of $H$ defined as, 
	\[\theta_i = \int t^{i} d H(t) \quad \textrm{for} \ i=1,2,3,4 .  \]	
	Define the first and second moments of the finite counterpart of limiting spectral distribution $S_T$ $F^{c_T, H_n}$ of the sample error covariance matrix $S_T$  as,
	\[\vartheta_1= \int x d F^{c_T, H_n},\quad \vartheta_2= \int x^2 d F^{c_T, H_n}. \]
	Suppose Assumptions 1, 2 and 4 hold for the fixed effects panel data model (1.1). Under the alternative $H_1$ and the LPA scheme, when$(n,T)\to \infty$, $c_T=n/T\to c\in(0,\infty)$, we have
	\begin{enumerate}
		\item If in addition $\gamma_4=3$, then
		\begin{eqnarray}
			s_1^{-1} \left[TU-T\left(\frac{\vartheta_2+\theta_2}{\vartheta_1^{2}}-1\right)\right] \overset{\mathcal{D}}{\longrightarrow} N(0,1),
		\end{eqnarray}
		where 
		\begin{eqnarray*}
			s_1^2 &\sim &  \frac{8\theta_4/c+4\theta_2^2+8c \theta_1^2 \theta_2+16\theta_1\theta_3}{\vartheta_1^4} -\frac{\left(16\theta_3/c+16\theta_1\theta_2\right)(\vartheta_2+\theta_2)}{\vartheta_1^5}\\
			&& \hspace{1.5cm}+\frac{8\theta_2(\vartheta_2+\theta_2)^2}{c \vartheta_1^6}.
		\end{eqnarray*}
		\item If $\Sigma_n$ is diagonal for all $n$ large, then
		\begin{eqnarray}
			s_2^{-1}\left[TU -T \left(\frac{\vartheta_2+\theta_2+\hat{\gamma}_4-3}{\vartheta_1^2}-1 \right) \right] \overset{\mathcal{D}}{\longrightarrow} N(0,1),
		\end{eqnarray}
		where
		\begin{eqnarray*}
			s_2^2 &\sim &  (\gamma_4-1)\Bigg( \frac{4\theta_4/c+2\theta_2^2+4 c \theta_1^2\theta_2+8\theta_1\theta_3}{\vartheta_1^4}+\frac{4\theta_2(\vartheta_2+\theta_2+\gamma_4-3)^2}{c \vartheta_1^6}  \\
			&& \hspace{3cm}-\frac{\left(8\theta_3/c+8\theta_1\theta_2\right)(\vartheta_2+\theta_2+\gamma_4-3)}{\vartheta_1^5}  \Bigg) .
		\end{eqnarray*}
	\end{enumerate}
\end{lemma}
When neither $\gamma_4=3$ nor $\Sigma_n$ is diagonal, the limiting distribution will become more complex, depending on the eigenvectors of $\Sigma_n$, see \cite{pan2008central}. From the above theorem, we know that the limiting distribution of John's test statistic $U$ depends on the moments of limiting spectral distribution of $\Sigma_n$ and the first two moments of $F^{c_T, H_n}$. 
\begin{corollary}\label{high_power_general}
	With the same assumptions as in Lemma \ref{high_alter}, when $(n,T)\to \infty, n/T=c_T\to c \in (0,\infty)$, the power of John's test under the alternative of general covariance $\Sigma_n$ is
	\item If in addition $\gamma_4=3$, then
	\begin{eqnarray*}
		P_1(H_1) \sim 1-\Phi \left(\frac{2 Z_\alpha +1+T\left(c+1 - \frac{\vartheta_2+\theta_2}{\vartheta_1^2} \right)}{s_1} \right),
	\end{eqnarray*}
	where $\alpha$ is the nominal test level, and $Z_\alpha, \Phi(\cdot)$ are the alpha upper quantile and cdf of standard normal distribution, respectively.
	\item If $\Sigma_n$ is diagonal for all $n$ large, then
	\begin{eqnarray*}
		P_2(H_1) \sim 1-\Phi \left(\frac{2 Z_\alpha +\hat{\gamma}_4-2+T\left(c+1 - \frac{\vartheta_2+\theta_2 +\gamma_4-3}{\vartheta_1^2} \right)}{s_2} \right).
	\end{eqnarray*}
\end{corollary}
Clearly, the consistency of John's test depends on the magnitude of $c+1$ compared to that of  $\frac{\vartheta_2+\theta_2+\gamma_4-3}{\vartheta_1^2}$. John's test is consistent when $c+1<\frac{\vartheta_2+\theta_2+\gamma_4-3}{\vartheta_1^2}$, otherwise inconsistent. 

\begin{conjecture}\label{conjection_general}
	With the same assumptions as in Lemma \ref{high_alter}, when $(n,T)\to \infty, n/T=c_T\to c \in (0,\infty)$, the difference between $\hat{U}$ and $U$ is
	\[T(\hat{U}-U )=O_p(1). \]
\end{conjecture}
It is natrual to get this conjecture since, as a special case of general bounded covariance matrix $\Sigma_n$, the difference $T(\hat{U}-U)$ is $O_p(1)$ under weak factor alternative.

Therefore, combining the above results, we can conclude that the consistency of the GRJ test is not guaranteed when the covariance matrix $\Sigma_n$ has a bounded norm under the LPA scheme.

\section{Proof of Lemma \ref{high_alter}}
Following the calculations in \cite{tian2015robust}, when $\gamma_4=3$, we get the following CLT.
\begin{proposition}
	Let $\{\ell_i \}_{1\leq i\leq n}$ be the eigenvalues of the sample error covariance matrix $S_T$, under the conditions of Lemma \ref{high_alter}, if $\gamma_4=3$, we have
	\begin{eqnarray*}
		T \left(\begin{array}{l}
			n^{-1}\sum_{i=1}^n \ell_i^2 -(\vartheta_2+\theta_2) \\ n^{-1}\sum_{i=1}^n \ell_i -\vartheta_1
		\end{array} \right) \overset{\mathcal{D}}{\longrightarrow} N\left( \left(\begin{array}{c} 0 \\0
		\end{array} \right), \mathbf{V_1} \right),
	\end{eqnarray*}
	with 
	\begin{eqnarray*}
		\mathbf{V_1} = \left(\begin{array}{cc}
			8 \theta_4/c+4\theta_2^2+8c \theta_1^2 \theta_2 +16\theta_1\theta_3 & 4\theta_3/c +4\theta_1\theta_2 \\ 4\theta_3/c +4\theta_1\theta_2 & 2\theta_2/c
		\end{array} \right).
	\end{eqnarray*}
\end{proposition}
Define the function $f(x,y) = \frac{x}{y^2}-1$, then $U= f\left(n^{-1} \sum_{i=1}^n \ell_i^2, n^{-1} \sum_{i=1}^n \ell_i \right)$. We have
\begin{eqnarray*}
	&&\frac{\partial f}{\partial x}\left(\vartheta_2+\theta_2, \vartheta_1 \right) = \frac{1}{\vartheta_1^2},\\ 
	&&\frac{\partial f}{\partial y}\left(\vartheta_2+\theta_2, \vartheta_1 \right) =  \frac{-2(\vartheta_2+\theta_2)}{\vartheta_1^3},\\
	&& f \left(\vartheta_2+\theta_2, \vartheta_1 \right) = \frac{\vartheta_2+\theta_2}{\vartheta_1^2}-1.
\end{eqnarray*}
By the delta method,
\begin{eqnarray*}
	T \left(U-f\left(\vartheta_2+\theta_2, \vartheta_1 \right) \right) \overset{\mathcal{D}}{\longrightarrow} N(0,s_1^2),
\end{eqnarray*}
where
\begin{eqnarray*}
	s_1^2 = \lim \left\{ \left( \begin{array}{c}
		\frac{\partial f}{\partial x}\left(\vartheta_2+\theta_2, \vartheta_1 \right) \\ \frac{\partial f}{\partial y}\left(\vartheta_2+\theta_2, \vartheta_1 \right)
	\end{array}
	\right)^T \cdot \mathbf{V}_1 \cdot \left( \begin{array}{c}
		\frac{\partial f}{\partial x}\left(\vartheta_2+\theta_2, \vartheta_1 \right) \\ \frac{\partial f}{\partial y}\left(\vartheta_2+\theta_2, \vartheta_1 \right)
	\end{array}
	\right) \right\}.
\end{eqnarray*}

For the second result, the central limit theorem in \cite{pan2008central} is applied to functions $f(x)=x^2$ and $g(x)=x$. 
\begin{proposition}
	Let $\{\ell_i \}_{1\leq i\leq n}$ be the eigenvalues of the sample covariance matrix $S_T$, under the conditions of Theorem \ref{high_alter}, if $\Sigma_n$ is diagonal for all $n$ large, we have
	\begin{eqnarray*}
		T \left(\begin{array}{l}
			n^{-1}\sum_{i=1}^n \ell_i^2 -(\vartheta_2+\theta_2+\gamma_4-3) \\ n^{-1}\sum_{i=1}^n \ell_i -\vartheta_1
		\end{array} \right) \overset{\mathcal{D}}{\longrightarrow} N\left( \left(\begin{array}{c} 0 \\0
		\end{array} \right), \mathbf{V_2} \right),
	\end{eqnarray*}
	with 
	\begin{eqnarray*}
		\mathbf{V_2} = (\gamma_4-1)\left(\begin{array}{cc}
			8 \theta_4/c+4\theta_2^2+8c \theta_1^2 \theta_2 +16\theta_1\theta_3 & 4\theta_3/c +4\theta_1\theta_2 \\ 4\theta_3/c +4\theta_1\theta_2 & 2\theta_2/c
		\end{array} \right).
	\end{eqnarray*}
\end{proposition}
Similarly, we can get,
\begin{eqnarray*}
	T \left(U-f\left(\vartheta_2+\theta_2+\gamma_4-3, \vartheta_1 \right) \right) \overset{\mathcal{D}}{\longrightarrow} N(0,s_2^2),
\end{eqnarray*}
where
\begin{eqnarray*}
	s_1^2 = \lim \left\{ \left( \begin{array}{c}
		\frac{\partial f}{\partial x}\left(\vartheta_2+\theta_2+\gamma_4-3, \vartheta_1 \right) \\ \frac{\partial f}{\partial y}\left(\vartheta_2+\theta_2+\gamma_4-3, \vartheta_1 \right)
	\end{array}
	\right)^T \cdot \mathbf{V}_2 \cdot \left( \begin{array}{c}
		\frac{\partial f}{\partial x}\left(\vartheta_2+\theta_2+\gamma_4-3, \vartheta_1 \right) \\ \frac{\partial f}{\partial y}\left(\vartheta_2+\theta_2+\gamma_4-3, \vartheta_1 \right)
	\end{array}
	\right) \right\}.
\end{eqnarray*}
The proof of Lemma \ref{high_alter} is complete. 

\bibliographystyle{plainnat}
\bibliography{reference}